\documentclass[12pt]{article}%
\usepackage{amsmath}
\usepackage{amsfonts}
\usepackage{amssymb}
\usepackage{graphicx}%
\graphicspath{{figures/}}
\usepackage{bbm}
\usepackage{amsthm}
\usepackage{enumerate}
\usepackage[T1]{fontenc}
\usepackage[utf8]{inputenc}
\usepackage{authblk}
\newtheorem{theorem}{Theorem}[section]
\newtheorem{corollary}{Corollary}[section]
\newtheorem{lemma}{Lemma}[section]
\RequirePackage[colorlinks,citecolor=blue,urlcolor=blue]{hyperref}

\newtheorem{remark}{Remark}[section]
\begin{document}
\numberwithin{equation}{section}
\title{The local backward heat problem}
\author{Thi Minh Nhat VO \thanks{Universit\'{e} d'Orl\'{e}ans, Laboratoire
MAPMO, CNRS UMR 7349, F\'{e}d\'{e}ration Denis Poisson, FR CNRS 2964,
B\^{a}timent de Math\'{e}matiques, B.P. 6759, 45067 Orl\'{e}ans Cedex 2,
France. Email address: vtmnhat@gmail.com.}}
\maketitle
\begin{abstract}
In this paper, we study the \textit{local backward problem} of a linear heat equation with time-dependent coefficients under the Dirichlet boundary condition. Precisely, we recover the initial data from the observation on a subdomain at some later time. Thanks to the `` optimal filtering '' method of Seidman, we can solve the \textit{ global backward problem}, which determines the solution at initial time from the known data on the whole domain. Then, by using a result of controllability at one point of time, we can connect local and global backward problem.\\

\end{abstract}
\bigskip
\section{Introduction and main result}
\subsection{Our motivation}
Inverse and ill-posed problems (see \cite{I}, \cite{P}, \cite{K}) are the heart of scientific inquiry and technological development. They play a significant role in engineering applications, as well as several practical areas, such as image processing, mathematical finance, physics, etc. and, more recently, modeling in the life sciences. During the last ten years or so, there have been remarkable developments both in the mathematical theory and applications of inverse problems. Especially, in various industrial purposes, for example steel industries, glass and polymer forming and nuclear power station, the "backward heat problem", which recovers the temperature in the heating system from the observation at some later time keeps an important position. On the other hand, from the mathematical point of view, it is well-known to be an ill-posed problem in the sense of Hadamard (see \cite{H}) due to the irreversibility of time. That is, there exists no solution from the given final data, and even if a solution exists, the small perturbations of the observation data may be dramatically scaled up in the solution. Hence, the interest of constructing some special regularization method is motivated. This topic has been studied extensively with lots of methods released such as: \textit{Tikhonov regularization} \cite{F}, \cite{M}, \cite{TS}, \cite{ZM}, \cite{MFH}, \textit{ Lavrentiev regularization} \cite{NT}, \cite{JSG}, \textit{truncation method} \cite{NTT}, \cite{KT}, \cite{ZFM}, \textit{filter method} \cite{S}, \cite{TKLT}, \cite{QW}, \textit{the quasi-boundary value method} \cite{DB}, \cite{KT}, \cite{QTTT} and other methods \cite{AE}, \cite{LL1}, \cite{LL2}, \cite{HX}, \cite{TQKT}, ... In \cite{S}, Seidman uses a so-called "optimal filtering" method in order to recover the solution at time $t >0$ with an optimal result. The idea of improving his outcome to reconstruct the solution at time $0$ is an interesting issue. Furthermore, the question that if we restrict our observation on a subregion inside the domain then how the local problem will be solved is also attractive. 
\subsection{Our problem}
 Let $\Omega$ be an open bounded domain in $\mathbb{R}^{n} (n\geq 1)$ with a boundary $\partial\Omega$ of class $C^{2}$; $T$ be a fixed positive constant. Let $p\in C^1\left([0;+\infty)\right)$ such that $0<p_1\leq p(t)\leq p_2, \forall t\in [0,+\infty)$, where $p_1$ and $p_2$ are some positive real numbers. Let $\omega$ be a nonempty, open subdomain of $\Omega$. We consider a linear heat equation with time-dependent coefficients, under the Dirichlet boundary condition with the state $u \in C^1\left((0,T);H_0^1(\Omega)\right)$:
\begin{equation}
\left\{\begin{array}{ll}
\partial_{t}u-p(t)\Delta u=0 & \text{in}~\Omega\times (0,T)\text{,}\\
u=0 & \text{on}~\partial \Omega \times (0,T)\text{.}\\
\end{array}\right.
\label{main_eq}
\end{equation}
Our target is constructing the initial solution $u(\cdot,0)$ when the local measurement data  of $u(\cdot,T)$ on the subdomain $\omega$ is available. In practice, the data at time $T$ is often measured by the physical instrument. Therefore, only a perturbed data $f$ can be obtained. Let $\delta > 0$ denote the noisy level with the following property
\begin{eqnarray}
\Vert u(\cdot,T)-f\Vert_{L^2(\omega)}\leq \delta \text{.}
\end{eqnarray}
Moreover, in order to assure the convergence of the regularization approximation to the initial data $u(\cdot,0)$, some a priori assumption on the exact solution is required
\begin{equation}
u(\cdot,0) \in H_0^1(\Omega)\text{.}
\label{u0}
\end{equation}
We will determine an approximate output $g$ of the unknown exact solution $u(\cdot,0)$ such that the error estimate $\mathfrak{e}(\delta)$ in
\begin{equation}
\Vert u(\cdot,0) - g \Vert_{L^2(\Omega)} \leq \mathfrak{e}(\delta)
\end{equation}
tends to $0$ when $\delta$ tends to $0$.
\subsection{ Relevant works}
Now, we consider how our problem can be solved so far in the past. In fact, there are lots of papers on the global backward problem but the works on the local one are restricted. There has been a sizeable literature on the special case $p \equiv 1$ with various methods. From now on, we will denote by $\delta$ the noisy level.
\begin{enumerate}
\item In 1996 (see \cite{S}), Seidman considers the heat equation which has the following form
\begin{equation}
\partial_t u - \nabla a\nabla u + q u = 0 ~~\text{on}~~ \Omega \times (0,T) ~~\text{with}~~u = 0 ~~\text{on}~~ \partial \Omega \times (0,T)
\end{equation}
where $a$ and $q$ belong to $L^\infty(\Omega)$. 
He succeeds in constructing the solution at a fixed time $t\in (0,T)$ from the observation $f$ satisfying $\Vert u(\cdot,T)-f\Vert_{L^2(\Omega)}\leq \delta$, under the assumption $u(\cdot,0) \in L^2(\Omega)$.  His strategy is using a "filter" with respect to the spectral decomposition of operator $\mathbf{A}: u \mapsto - \nabla a\nabla u + q u$, which is defined as
\begin{equation}
\mathbf{F}(t)e_i = \min\left\{1, e^{-\lambda_i (T-t)} \frac tT\left(\frac{\Vert u(\cdot,0) \Vert_{L^2(\Omega)}}{\delta}\right)^{1-\frac tT}\right\}e_i
\end{equation}
where $\{\lambda_i\}_{i\geq1}$ and $\{e_i\}_{i\geq1}$ are respectively denoted by the eigenvalues and the corresponding eigenfunctions of the operator $\mathbf{A}$. Then, 
he can get the optimal result, which is
\begin{equation} \Vert u(\cdot,t) - g_t \Vert_{L^2(\Omega)} \leq \delta^{\frac tT}\Vert u(\cdot,0) \Vert_{L^2(\Omega)}^{1-\frac tT} \text{.}
\end{equation}
The regularization solution $g_t$ at time $t$ is constructed as
\begin{equation}
g_t := \sum \limits_{i=1}^\infty e^{\lambda_i(T-t)}\left(\int_\Omega f(x)e_i(x)dx\right)\mathbf{F}(t)e_i \text{.}
\end{equation}
\item By generalizing the result of Seidman, in \cite{TS}, Tautenhahn and Schr\"{o}ter provide us a definition of the term "optimal method", in a sense the error of the estimate between the exact solution and the approximate one defined from the optimal method can not be greater than the best possible worst case error (see Definition 1.1, page 478).  Their interest is finding the optimal results in different regularization methods for solving the backward heat equations. According to this sense, the result of Seidman is optimal.
\item In 2007, Trong et al. (see \cite{TQKT}) improve the quasi-boundary value method to regularize the 1D backward heat equation. They succeed in recovering the initial data with the following error:
\begin{equation}\left\Vert u(\cdot,0) - g \right\Vert _{L^2(\Omega)}  \leq \sqrt[4]{8}C\sqrt[4]{T}\left(\ln \frac{\Vert u(\cdot,0)\Vert_{L^2(\Omega)}}\delta\right)^{-\frac 14}\end{equation}
where $C$ is a positive constant depending on $\Vert u(\cdot,0) \Vert_{H_0^1(\Omega)}$ (see also \cite{AP}).
 \end{enumerate}
 The problem with case $p \not\equiv constant$ is recently concerned, which can be mentioned in some following writings:
\begin{enumerate}
\item In 2013, Tuan et al. in \cite{TQTT} consider the 1D backward heat equation with time-dependent coefficients. They use a so-called "modified method" to get the result below
\begin{equation} \Vert u(\cdot,t) - g_t \Vert_{L^2(\Omega)} \leq \left(1+\Vert u(\cdot,0)\Vert_{L^2(\Omega)}\right)\left(\frac{\delta}{\Vert u(\cdot,0)\Vert_{L^2(\Omega)}}\right)^{\frac{p_1^2 t}{p_2^2 T}}\text{.}\end{equation}
\item In 2014, Zhang, Fu and Ma (see  \cite{ZFM}) also study on the 1D backward heat equation with time-dependent coefficients, but use the truncation method. They can recover the solution at time $t \in \left(T(1-\frac{p_1}{p_2});T\right)$ satisfying
\begin{equation} \left\Vert u(\cdot,t) - g_t \right\Vert _{L^2(\Omega)}  \leq \Vert u(\cdot,0) \Vert_{L^2(\Omega)}^{1-\frac tT}((\tau +1)\delta)^{\frac tT}+\left( \frac{ \Vert u(\cdot,0) \Vert_{L^2(\Omega)}}{\tau -1}\right)^{\frac{p_2(T-t)}{p_1T}}\delta^{\frac{(p_2-p_1)T+p_2t}{p_1T}}\end{equation}
for some constant $\tau >1$.
\item In 2016,  Khanh and Tuan (in \cite{KT}) solve an initial inverse problem for an inhomogeneous heat equation by using high frequency truncation method. Under the assumption that $u(\cdot,0)\in H_0^1(\Omega)$, they can recover the initial data with the below error:
\begin{equation}
\left\Vert u(\cdot,0) - g \right\Vert _{L^2(\Omega)}  \leq \left(\frac{\delta}{\Vert u(\cdot,0)\Vert_{L^2(\Omega)}}\right)^{\frac 1{2T}}\frac{\sqrt{\ln\left(\frac{\Vert u(\cdot,0)\Vert_{L^2(\Omega)}}{\delta}\right)}}{\sqrt{2p_2T}}+\frac{\sqrt{2p_2T}\Vert u(\cdot,0)\Vert_{H_0^1(\Omega)}}{\sqrt{\ln\left(\frac{\Vert u(\cdot,0)\Vert_{L^2(\Omega)}}{\delta}\right)}}\text{.}
\end{equation}
\end{enumerate}
For the local inverse problem, we can pick up some of following works:
\begin{enumerate}
\item In 1995, Yamamoto in \cite{Y} proposes a reconstruction formula for the spatial dependence of the source term in a wave equation $\partial_{tt}u-\Delta u = f(x)\sigma(t)$, assuming $\sigma(t)$ known, from local measurement using exact controllability. 
\item In 2009, Li, Yamamoto and Zou in \cite{LYZ} study the conditional stability of inverse problems. Here, the known data is observed in a subregion along a time period which may  start at some point, possibly far away from the initial data.
\item In 2011, Garc\'ia and Takahashi (see \cite{GT}) present some abstract results of a general connection between null-controllability and several inverse problems for a class of parabolic equations.
\item In 2013, Garc\'ia, Osses and Tapia in \cite{GOT} succeed in determining the heat source from a single internal measurement of the solution, thanks to a family of null control functions.
\end{enumerate}
\subsection{Our method of solving the global backward problem (GBP) and the local backward problem (LBP)}
Firstly, we deal with the global backward problem (GBP), which recovers the initial data from the observation on the whole domain at some later time $\tau>0$. Here, we assume that there exists solution of  the (\ref{main_eq}) satisfying the a priori condition (\ref{u0}) and $\bar{f}$ be the known data on $\Omega$ at time $\tau$ such that $\Vert u(\cdot,\tau)-\bar{f}\Vert_{L^2(\Omega)}\leq \delta$ for some $\delta >0$. We will determine a function $g$ which approximate the initial data. Our idea of constructing such function $g$ is from the ``optimal filtering method'' of Seidman (see \cite{S}): First, we define a continuous operator depending on a regularization parameter $\alpha$
\begin{eqnarray}
\mathcal{R}_{\alpha} : L^2(\Omega) & \to & L^2(\Omega) \notag \\
    \phi & \mapsto &  \sum\limits_{i=1}^\infty \min\{e^{\lambda_i\int_0^\tau p(s)ds};\alpha\}\left(\int_\Omega \phi(x)e_i(x)dx\right)e_i \text{;}
    \notag
 \end{eqnarray} 
Then, the function $\mathcal{R}_{\alpha}\bar{f}$ will be closed to the exact solution $u(\cdot,0)$ in $L^2(\Omega)$ where $\alpha$ is the minimizer of the problem $\min\limits_{\alpha>0}\Vert u(\cdot,0) - \mathcal{R}_\alpha \bar{f}\Vert_{L^2(\Omega)}$.\\
Secondly,  for the local backward problem (LBP), whose observation is measured on a subdomain, we need to use a tool of controllability to link with the (GBP). Precisely, we use the assertion about the existence of a sequence of control functions to get the information of solution on the whole domain $\Omega$ from the given data on the subdomain $\omega$: For each $i=1,2,...$, for any $\varepsilon >0$, there exists $h_i\in L^2(\omega)$ such that the solution of 
\begin{equation}
\left\{\begin{array}{ll}
\partial_{t}\varphi _i-p(t)\Delta \varphi _i=0 & \text{in}~\Omega\times (0,2T)\setminus\{T\},\\
\varphi _i=0 & \text{on}~\partial \Omega \times (0,2T),\\
\varphi _i(\cdot,0)=e_i& \text{in}~ \Omega \text{,}\\
\varphi _i(\cdot,T)=\varphi _i(\cdot,T^{-}) + \mathbbm{1}_{\omega}{h_i} &\text{in} ~\Omega
\end{array}\right.
\label{vi}
\end{equation}
satisfies $\left\Vert \varphi _i(\cdot,2T) \right\Vert_{L^2(\Omega)}\leq \varepsilon$. Here, $\mathbbm{1}_\omega$ presents for the characteristic function on the region $\omega$ and $\varphi_i(T^-)$ denotes the left limit of the function $\varphi_i$ at time $T$. Multiplying $\partial_{t}\varphi _i-p(t)\Delta \varphi _i=0$ by $u(\cdot, 2T-t)$ and using some computation technique, we can get the approximate solution $\bar{f}$ at time $\tau = 3T$, which is \[\Vert u(\cdot, 3T) - \bar{f}\Vert_{L^2(\Omega)}\leq \mathcal{E}(\delta)\text{.}\] Here, $\bar{f}$ is computed by known data $h_i$ and $f$ and $\mathcal{E}(\delta)$ is a function of $\delta$ such that $\mathcal{E}(\delta)\rightarrow 0$ when $\delta \rightarrow 0$. Lastly, applying the result of (GBP) with the information at $3T$ on the whole domain, the initial data of (\ref{main_eq}) is reconstructed.
\subsection{Spectral theory}
As a direct consequence of spectral theorem for compact, self-adjoint operators (see Theorem 9.16, page 225, \cite{HN}), there exists a sequence of positive real eigenvalues of the operator $-\Delta$, which denoted by $\{\lambda_i\}_{i=1,2,...}$ where 
\begin{equation}
\left\{\begin{array}{ll}
0<\lambda_1\leq\lambda_2\leq\lambda_3\leq....\text{,}\\
\lambda_i \rightarrow \infty ~\text{as}~ i \rightarrow \infty \text{.}
\end{array}\right.
\end{equation}
Moreover, there exists an orthonormal basis $\{e_i\}_{i=1,2,...}$ of $L^2(\Omega)$, where $e_i \in H_0^1(\Omega)$ is an eigenfunction corresponding to $\lambda_i$
\begin{equation}
\left\{\begin{array}{ll}
-\Delta e_i = \lambda_i e_i &\text{in}~ \Omega\text{,}\\
e_i = 0 &\text{on}~ \partial \Omega \text{.}
\end{array}\right.
\end{equation}
When $u_0 \in L^2(\Omega)$ and $u_0 = \sum\limits_{i=1}^{\infty} a_ie_i$ with $a_i=\int_\Omega u_0(x)e_i(x)dx$ and $\sum\limits_{i=1}^{\infty} \left|a_i\right|^2<\infty$, then
\begin{equation}
u(\cdot,t) = \sum\limits_{i=1}^{\infty} a_ie^{-\lambda_i\int_0^t p(s)ds}e_i
\end{equation}
is the unique solution of
\begin{equation}
\left\{\begin{array}{ll}

\partial_{t}u-p(t)\Delta u=0 & \text{in}~\Omega\times (0,T),\\
u=0 & \text{on}~\partial \Omega \times (0,T),\\
u(\cdot,0)=u_0&\text{in}~ \Omega.
	
\end{array}\right.
\end{equation}
\subsection{Main result}
\begin {theorem}
Let $u$ be the solution of (\ref{main_eq}) with the a priori bound (\ref{u0}). Let $f\in L^2(\omega)$ and  $0<\delta<\Vert u(\cdot,0)\Vert_{L^2(\Omega)}$ satisfying
\begin{eqnarray}
\Vert u(\cdot,T)-f\Vert_{L^2(\omega)}\leq \delta \text{.}
\label{delta1}
\end{eqnarray}
There exists a function $g \in L^2(\Omega)$ and a constant $C=C(\Omega, \omega, p)>0$ such that the following estimate holds
\begin{equation}
\left\Vert u(\cdot,0) - g \right\Vert _{L^2(\Omega)}  \leq \frac{Ce^{\frac CT}\sqrt{T}\Vert u(\cdot, 0)\Vert_{H_0^1(\Omega)}}{\sqrt{\ln \frac{\Vert u(\cdot,0)\Vert_{L^2(\Omega)}}{\delta}}}\text{.}\label{maines}
\end{equation}
\label{local}
\end{theorem}
\begin{remark}
\begin{enumerate}
\item When $\delta < De^{-D\left(T+\frac 1T\right)}\Vert u(\cdot,0)\Vert_{L^2(\Omega)}$ for some positive constant $D=D(\Omega, \omega, p)$, the approximate solution of the initial data satisfying (\ref{maines}) is constructed as below
\begin{eqnarray}
g&:= &- \sum\limits_{i=1}^\infty \min\{e^{\lambda_i\int_0^{3T}p(s)ds},\alpha\}e^{-\lambda_i\int_{2T}^{3T}p(s)ds}\left(\int_\omega h_i(x)f(x)dx\right)e_i
\end{eqnarray}
where $\{h_i\}_{i\geq 1}$ is a sequence of control functions  (see Section 4) and $\alpha$ is the regularization parameter given by
\begin{equation}
\alpha = \mathcal{A}\left(\mathcal{B}^{-1}\left(\frac{\sqrt{3p_2T}\Vert u(\cdot,0) \Vert_{H_0^1(\Omega)}}{K_1e^{\frac {K_1}T}{\Vert u(\cdot,0)\Vert_{L^2(\Omega)}}^{1-k_1}\delta^{k_1}}\right)\right)
\end{equation}
with 

\begin{enumerate}[(i)]
  \item  \begin{eqnarray}
\mathcal{A} : \left[0;+\infty\right) & \to & \left[\frac{\sqrt{e}}{2}; +\infty\right) \notag \\
    x & \mapsto &  \frac{e^x}{1+2x} \text{,}
    \label{funcI}
 \end{eqnarray} 
 \item \begin{eqnarray}
\mathcal{B} : (0;+\infty) & \to & (0; +\infty) \notag \\
    x & \mapsto &  \sqrt{x}e^x
    \label{funcH}
 \end{eqnarray}
  The existence of the function $\mathcal{B}^{-1}$ dues to the bijection property of the function $\mathcal{B}$ on $(0,+\infty)$,
\item $K_1=K_1(\Omega, \omega, p)>1$ and $k_1=k_1(\Omega, \omega, p)\in (0,1)$. All these constants can be explicitly computed when $\Omega$ is convex or star-shaped with respect to some $x_0 \in \Omega$.
\end{enumerate}
\item The estimate (\ref{maines}) connects to the well-known following estimate
\begin{equation}
\left\Vert u(\cdot,0)\right\Vert _{L^2(\Omega)}  \leq \frac{C\sqrt{1+T+\frac{1}{T}}\Vert u(\cdot, 0)\Vert_{H_0^1(\Omega)}}{\sqrt{\ln \frac{\Vert u(\cdot,0)\Vert_{L^2(\Omega)}}{\Vert u(\cdot,T)\Vert_{L^2(\omega)}}}}\text{.}
\label{appendix}
\end{equation}
for some positive constant $C=C(\Omega,\omega,p)$ (see Appendix for the proof).
\end{enumerate}
\label{remark2}
\end{remark}
\subsection{Outline}
Section 2 will give us a result of the (GBP) (see Theorem \ref{global}), where the known data is observed on the whole domain. In section 3, we construct an observation estimate at one point of time for the parabolic equations with time-dependent coefficients (see Theorem \ref{nonconvex} and Theorem \ref{convex}). This is an important preliminary of the approximate controllability (see Theorem \ref{control}), which is studied in section 4. Lastly, combining the result of controllability and global backward, we get the proof of the Theorem \ref{local}, mentioned in section 5.

\section{Global backward problem}
First of all, we need to consider the special case, that is $\omega \equiv \Omega$. In \cite{S}, Seidman succeeds in recovering the solution at time $t >0$ by an optimal filtering method, under the a priori condition $u(\cdot,0)\in L^2(\Omega)$. Here we will use his method to recover the initial data at time $0$ but with the stronger assumption $u(\cdot,0)\in H_0^1(\Omega)$.
\begin {theorem}
Let $u$ be the solution of (\ref{main_eq}) satisfying the a priori condition (\ref{u0}). Let $\bar{f}\in L^2(\Omega)$ and $\delta >0$ having the following property
 \begin{equation}
\Vert u(\cdot,T)-\bar{f}\Vert_{L^2(\Omega)}\leq \delta \text{.}
\label{globalz}
\end{equation}
There exists a function $g \in L^2(\Omega)$ such that for any $\zeta > \frac{\delta^2}{2\lambda_1p_2T\Vert u(\cdot,0)\Vert^2_{L^2(\Omega)}}$, the following estimate holds
\begin{equation}
\left\Vert u(\cdot,0) - g \right\Vert _{L^2(\Omega)}  \leq \frac{\sqrt{(1+\zeta) p_2T}\left\Vert u(\cdot,0)\right\Vert_{H_0^1(\Omega)}}{\sqrt{\ln \left(\sqrt{2\zeta\lambda_1 p_2T}\frac{\left\Vert u(\cdot,0)\right\Vert_{L^2(\Omega)}}{\delta}\right)}} \text{.}
\label{ln}
\end{equation}
\label{global}
\end{theorem}
\begin{remark}
\begin{enumerate}
\item When $\delta < \Vert u(\cdot,0)\Vert_{L^2(\Omega)}e^{-\lambda_1p_2T}$, the approximate solution of the initial data satisfying (\ref{ln}) is constructed as
\begin{equation} g := \sum\limits_{i=1}^{\infty} \min\left\{e^{\lambda_i\int_0^T p(s)ds},\bar{\alpha}\right\}\int_\Omega \bar{f}(x)e_i(x)dxe_i \text{.}\end{equation}
Here, the regularization parameter $\bar{\alpha}$ is given by 
\begin{equation}
\bar{\alpha}:=\mathcal{A}\left(\mathcal{B}^{-1}\left(\frac{\sqrt{p_2T}\left\Vert u(\cdot,0)\right\Vert_{H_0^1(\Omega)}}{\delta}\right)\right)
\end{equation}
with $\mathcal{A}$ and $\mathcal{B}$ being respectively defined in (\ref{funcI}) and (\ref{funcH}).
\item When $\delta < \Vert u(\cdot,0)\Vert_{L^2(\Omega)}$, we can choose $\zeta = \frac{1}{2\lambda_1p_2T}$ in order to get
\begin{equation}
\left\Vert u(\cdot,0) - g \right\Vert _{L^2(\Omega)}  \leq \frac{\sqrt{p_2T+\frac{1}{2\lambda_1}} \left\Vert u(\cdot,0)\right\Vert_{H_0^1(\Omega)}}{\sqrt{\ln \frac{\left\Vert u(\cdot,0)\right\Vert_{L^2(\Omega)}}{\delta}}} \text{.}
\end{equation}
This connects to the well-known following estimate
\begin{equation}
\left\Vert u(\cdot,0)\right\Vert _{L^2(\Omega)}  \leq \frac{\sqrt{p_2T}\left\Vert u(\cdot,0)\right\Vert_{H_0^1(\Omega)}}{\sqrt{\ln \frac{\left\Vert u(\cdot,0)\right\Vert_{L^2(\Omega)}}{\Vert u(\cdot,T)\Vert_{L^2(\Omega)}}}}\text{.}
\end{equation}
\end{enumerate}
\label{remark_global}
\end{remark}
\textbf{Proof of Theorem \ref{global}}
\begin{proof}
For the case $\delta \geq \Vert u(\cdot,0)\Vert_{L^2(\Omega)}e^{-\lambda_1p_2T}$, the estimate (\ref{ln}) holds with $g=0$. Indeed, combining with the fact that $\sqrt{2\zeta\lambda_1p_2T}\leq e^{\zeta\lambda_1p_2T}~~\forall \zeta >0$, we get 
\begin{equation} 
\sqrt{2\zeta\lambda_1p_2T}\frac{\Vert u(\cdot,0)\Vert_{L^2(\Omega)}}{\delta}\leq e^{(1+\zeta)\lambda_1p_2T}
\text{.}\end{equation}
It implies that
\[\frac 1{\sqrt{\lambda_1}} \leq \frac {\sqrt{(1+\zeta)p_2T}}{\sqrt{\ln \sqrt{2\zeta\lambda_1p_2T}\frac{\Vert u(\cdot,0)\Vert_{L^2(\Omega)}}{\delta}}}\text{.}\]
Hence
\begin{eqnarray}
\left\Vert u(\cdot,0)\right\Vert _{L^2(\Omega)} &\leq&  \frac{\Vert u(\cdot,0)\Vert_{H_0^1(\Omega)}}{\sqrt{\lambda_1}}\notag\\
&\leq& \frac {\sqrt{(1+\zeta)p_2T}\Vert u(\cdot,0)\Vert_{H_0^1(\Omega)}}{\sqrt{\ln \sqrt{2\lambda_1p_2T}\frac{\Vert u(\cdot,0)\Vert_{L^2(\Omega)}}{\delta}}}\text{.}
\end{eqnarray}
The main purpose concerns the case when 
\begin{equation}
\delta < \Vert u(\cdot,0)\Vert_{L^2(\Omega)}e^{-\lambda_1p_2T}\text{.}
\label{delta}
\end{equation}
In this case, we will determine the regularization solution at time $0$ as follows: First of all, Step 1 will provide us the construction of a continuous operator $\mathcal{R}_{\beta}$ depending on a parameter $\beta$, which will be chosen later. The regularization solution $g$ is defined by applying this operator on the known-data $\bar{f}$. Secondly, in Step 2, we compute the error between the exact solution and the approximate solution defined in Step 1. Lastly, by minimizing the error in Step 2 with respect to $\beta$, we can obtain the final result in Step 3.

\textbf{Step 1:} Construct the regularization solution. \\
Let us define a continuous function $\mathcal{R}_{\beta}$ depending on a positive parameter $\beta$, which will be chosen later:
\begin{eqnarray}
\mathcal{R}_{\beta} : L^2(\Omega) & \to & L^2(\Omega) \notag \\
    f & \mapsto &  \sum\limits_{i=1}^\infty \min\{e^{\lambda_i\int_0^Tp(s)ds};\beta\}\left(\int_\Omega f(x)e_i(x)dx\right)e_i
    \notag
 \end{eqnarray} 
Put $g:=\mathcal{R}_{\beta} \bar{f}$. We will prove that such defined function $g$ approximate the exact solution $u(\cdot,0)$ with some suitable choice of $\beta$.\\

\textbf{Step 2:} Compute the error $\Vert u(\cdot,0) - g \Vert_{L^2(\Omega)}$.\\
Put $g_T:=\mathcal{R}_{\beta} u(\cdot,T)$, we will compute the error by using the following triangle inequality
\begin{eqnarray}
\Vert u(\cdot,0) - g \Vert_{L^2(\Omega)} \leq \Vert u(\cdot,0) - g_T \Vert_{L^2(\Omega)} + \Vert g_T - g \Vert_{L^2(\Omega)}\text{.}
\end{eqnarray}
On one hand, we have
\begin{eqnarray}
\left\Vert g - g_T \right\Vert_{L^2(\Omega)}
&=& \left\Vert\sum\limits_{i=1}^{\infty} \min\{e^{\lambda_i\int_0^Tp(s)ds}; \beta\}\int_\Omega \left(\bar{f}(x)-u(x,T)\right)dx e_i\right\Vert_{L^2(\Omega)}\notag\\
&\leq& \beta\left\Vert \bar{f}-u(\cdot,T)\right\Vert_{L^2(\Omega)}\leq \beta\delta \text{.}
\end{eqnarray}
On the other hand 
\begin{eqnarray}
&~&\left\Vert u(\cdot,0) - g_T \right\Vert_{L^2(\Omega)} \notag\\
&=& \left\Vert \sum\limits_{i=1}^{\infty} \int_{\Omega} u(x,0)e_i(x)dx e_i - \sum_{i =1}^{\infty} \min\{e^{\lambda_i\int_0^Tp(s)ds}; \beta\}e^{-\lambda_i\int_0^Tp(s)ds}\int_\Omega u(x,0)e_i(x)dx e_i \right\Vert_{L^2(\Omega)}\notag\\
&=&  \left\Vert \sum\limits_{i=1}^{\infty} \left(1- \min\{e^{\lambda_i\int_0^Tp(s)ds}; \beta\}e^{-\lambda_i\int_0^Tp(s)ds}\right)\int_\Omega u(x,0)e_i(x)dx e_i \right\Vert_{L^2(\Omega)}\notag\\
&=&  \left\Vert \sum_{e^{\lambda_i\int_0^Tp(s)ds}>\beta} \left(1- \beta e^{-\lambda_i\int_0^Tp(s)ds}\right)\int_\Omega u(x,0)e_i(x)dx e_i \right\Vert_{L^2(\Omega)}\notag\\
&=& \left\Vert \sum_{e^{\lambda_i\int_0^Tp(s)ds}>\beta} \frac{\left(1- \beta e^{-\lambda_i\int_0^Tp(s)ds}\right)}{\sqrt{\lambda_i}}\sqrt{\lambda_i}\int_\Omega u(x,0)e_i(x)dx e_i \right\Vert_{L^2(\Omega)}\notag\\
&\leq& \sup_{\lambda\geq\lambda_1} \frac{\left(1 - \beta e^{-\lambda\int_0^Tp(s)ds}\right)}{\sqrt{\lambda}}\left\Vert u(\cdot,0)\right\Vert _{H_0^1(\Omega)} \notag\\
&\leq& \sup_{\lambda\geq\lambda_1} \frac{(1 - \beta e^{-\lambda p_2T} )}{\sqrt{\lambda p_2T}}\sqrt{p_2T}\left\Vert u(\cdot,0)\right\Vert _{H_0^1(\Omega)} \text{.}
\end{eqnarray}
Now, we solve the problem of finding $\sup\limits_{\lambda\geq\lambda_1} \frac{(1 - \beta e^{-\lambda p_2T} )}{\sqrt{\lambda p_2T}}$.\\
Define \begin{eqnarray}
\mathcal{F} : [\lambda_1;+\infty) & \to & (0; +\infty) \notag \\
    \lambda & \mapsto &  \frac{(1 - \beta e^{-\lambda p_2T} )}{\sqrt{\lambda p_2T}}\text{.}
    \notag
 \end{eqnarray} 
Obviously, $\mathcal{F}$ is differentiable and 
\begin{equation}
\mathcal{F}'(\lambda) = \frac{\beta p_2T e^{-\lambda p_2T}(1+2\lambda p_2T)-p_2T}{2\sqrt{\lambda p_2T}\lambda p_2T} \text{.}\end{equation}
The equation $\mathcal{F}'(\lambda)=0$ is equivalent to
 \begin{equation}
 \beta = \frac{e^{\lambda p_2T}}{1+2\lambda p_2T}\text{.}
 \label{beta1}
 \end{equation}
We will choose $\beta$ such that the equation (\ref{beta1}) has a unique solution $\bar{\lambda} \geq \lambda_1$.
Let us remind the function $\mathcal{A}$ defined in (\ref{funcI}):
 \begin{eqnarray}
\mathcal{A} : \left[0;+\infty\right) & \to & \left[\frac{\sqrt{e}}{2}; +\infty\right) \notag \\
    x & \mapsto &  \frac{e^x}{1+2x} \text{.}
    \notag
 \end{eqnarray} 
Note that the equation (\ref{beta1}) has a unique solution $\bar{\lambda} \geq \lambda_1$ if and only if 
\begin{equation}
\beta > \mathcal{A}\left(\lambda_1 p_2T\right) \text{.}
\label{beta2}
\end{equation}
Suppose the condition (\ref{beta2}) is satisfied then there exists a unique $\bar{\lambda}\geq \lambda_1$ such that $\mathcal{F}'(\bar{\lambda})=0$ and $\beta = \mathcal{A}(\bar{\lambda}p_2T)$. We can write $\bar{\lambda}p_2T=\mathcal{A}^{-1}(\beta)$. On the other hand, the fact that $\mathcal{F}'(\lambda_1)>0$ leads us to the conclusion: the function $\mathcal{F}$ is strictly increasing on $(\lambda_1,\bar{\lambda})$ and strictly decreasing on $(\bar{\lambda}, +\infty)$. Consequently, $\mathcal{F}$ gets supremum at $\bar{\lambda}$, i.e
\begin{equation}
 \mathcal{F}(\bar{\lambda}) = \sup\limits_{\lambda\geq\lambda_1} \mathcal{F}(\lambda)\text{.}
 \end{equation}
 
\textbf{Step 3:} Minimize the error with a suitable choice of $\beta$.\\
Combining the two above steps, we get
\begin{eqnarray}
\left\Vert u(\cdot,0) - g \right\Vert _{L^2(\Omega)} &\leq &\beta \delta + \frac{(1 - \beta e^{-\bar{\lambda} p_2T} )}{\sqrt{\bar{\lambda} p_2T}} \sqrt{p_2T}\left\Vert u(\cdot,0)\right\Vert _{H_0^1(\Omega)}\notag\\
&=& \Theta \delta e^{\mathcal{A}^{-1}(\beta)} + (1-\Theta)\frac{\sqrt{p_2T}\left\Vert u(\cdot,0) \right\Vert_{H_0^1(\Omega)}}{\sqrt{\mathcal{A}^{-1}(\beta)}}
\label{fou}
\end{eqnarray}
where $\Theta = \beta e^{-\bar{\lambda} p_2T }$. Note that 
\[
\Theta \mathsf{A} +(1-\Theta) \mathsf{B} \geq \min\{\mathsf{A}, \mathsf{B}\}~~ \forall ~~\mathsf{A},\mathsf{B} >0, \Theta \in (0,1)\text{.}\]
The equality occurs when and only when $\mathsf{A}=\mathsf{B}$. Hence, in order to minimize the right-hand side of (\ref{fou}), we will choose $\beta$ such that 
\begin{equation}
  \delta e^{\mathcal{A}^{-1}(\beta)}= \frac{\sqrt{p_2T}\left\Vert u(\cdot,0) \right\Vert_{H_0^1(\Omega)}}{\sqrt{\mathcal{A}^{-1}(\beta)}} \text{.}
\label{beta}
\end{equation}
The choice of $\beta =\bar{\alpha} := \mathcal{A}\left(\mathcal{B}^{-1} \left(\frac{\sqrt{p_2T}\left\Vert u(\cdot,0)\right\Vert_{H_0^1(\Omega)}}{\delta}\right)\right)$ satisfies the condition (\ref{beta2}) (dues to the assumption (\ref{delta}) on the smallness of $\delta$) and the estimate (\ref{beta}). Therefore, we get the following estimate
\begin{eqnarray}
\left\Vert u(\cdot,0) - g \right\Vert _{L^2(\Omega)} 
&\leq& \frac{\sqrt{p_2T}\left\Vert u(\cdot,0) \right\Vert_{H_0^1(\Omega)}}{\sqrt{\mathcal{A}^{-1}(\bar{\alpha})}}\notag\\
&\leq&\frac{\sqrt{p_2T}\left\Vert u(\cdot,0) \right\Vert_{H_0^1(\Omega)}}{\sqrt{\mathcal{B}^{-1} \left(\frac{\sqrt{p_2T}\left\Vert u(\cdot,0)\right\Vert_{H_0^1(\Omega)}}{\delta}\right)}} \text{.}
\label{up}
\end{eqnarray}
Due to the definition of the function $\mathcal{B}$ (see \ref{funcH}), (\ref{up}) becomes
\begin{eqnarray}
\frac{\sqrt{p_2T}\left\Vert u(\cdot,0)\right\Vert_{H_0^1(\Omega)}}{\delta} \leq \frac{\sqrt{p_2T}\left\Vert u(\cdot,0)\right\Vert_{H_0^1(\Omega)}}{\left\Vert u(\cdot,0) - g\right\Vert_{L^2(\Omega)}}e^{\frac{p_2T \left\Vert u(\cdot,0)\right\Vert^2_{H_0^1(\Omega)}}{\left\Vert u(\cdot,0) - g\right\Vert^2_{L^2(\Omega)}}} \text{.}
\end{eqnarray}
Using the fact that $\sqrt{2\zeta}x\leq e^{\zeta x^2} ~~\forall \zeta >0~~\forall x>0$, one obtains
\begin{eqnarray}
\frac{\sqrt{p_2T}\left\Vert u(\cdot,0)\right\Vert_{H_0^1(\Omega)}}{\delta} \leq \frac 1{\sqrt{2\zeta}}e^{\frac{(1+\zeta)p_2T \left\Vert u(\cdot,0)\right\Vert^2_{H_0^1(\Omega)}}{\left\Vert u(\cdot,0) - g\right\Vert^2_{L^2(\Omega)}}} \text{.}
\end{eqnarray}
It is equivalent to
\begin{equation}
\frac{(1+\zeta)p_2T \left\Vert u(\cdot,0)\right\Vert^2_{H_0^1(\Omega)}}{\left\Vert u(\cdot,0) - g\right\Vert^2_{L^2(\Omega)}} \geq \ln \left(\frac{\sqrt{2\zeta p_2T}\left\Vert u(\cdot,0)\right\Vert_{H_0^1(\Omega)}}{\delta}\right)\text{.}
\end{equation}
With $\zeta > \frac{\delta^2}{2\lambda_1p_2T\Vert u(\cdot,0)\Vert^2_{L^2(\Omega)}}$, it deduces that
\begin{equation}
\left\Vert u(\cdot,0) - g \right\Vert _{L^2(\Omega)}  \leq \frac{\sqrt{(1+\zeta)p_2T}\left\Vert u(\cdot,0)\right\Vert_{H_0^1(\Omega)}}{\sqrt{\ln \left(\sqrt{2\zeta\lambda_1 p_2T}\frac{\left\Vert u(\cdot,0)\right\Vert_{L^2(\Omega)}}{\delta}\right)}} \text{.}
\end{equation}
\end{proof}
For the local case $\omega \Subset \Omega$, it is required the existence of control functions on the subdomain at some point of time in order to link with a global result. This controllability problem has a sustainable connection with the observability one, which will be studied in the next Section.  
 \section{ Observability at one point of time}
The issue on constructing an observation estimate is widely studied. It can be solved by global Carleman inequality, which is presented in \cite{FI}; by using the estimate of Lebeau and Robbiano (see \cite{LR}) or by transmutation (see \cite{EZ}). Recently, Phung et al. provide a different method which is based on properties of the heat kernel with a parametric of order 0. In \cite{PW1} and \cite{PW2}, the authors work on a linear equation which has form
\begin{equation*}
\partial_t v - \Delta v + av + b\nabla v = 0  ~~\text{in} ~~\Omega \times (0,T)\text{.}
\end{equation*}
Here, $a\in L^\infty((0,T),L^q(\Omega))$ with $q\geq 2$ if $n=1$ and $q>n$ if $n\geq 2$; $b\in L^\infty(\Omega \times (0,T))^n$ and $\Omega$ must be convex. Then, by using some geometrical techniques, Phung et. al improve their previous results by working on a general domain (i.e $\Omega$ is even convex or not). For the following form of linear equation
\begin{equation*}
\partial_t v - \Delta v + av = 0  ~~\text{in} ~~\Omega \times (0,T)
\end{equation*}
where $a\in L^\infty(\Omega\times(0,T))$, see \cite{PWZ}. For the parabolic equations with space-time coefficients
\begin{equation*}
\partial_t v - \nabla(A\nabla v) + av + b\nabla v= 0  ~~\text{in} ~~\Omega \times (0,T)
\end{equation*}
where $a\in L^\infty(\Omega\times(0,T))$,  $b\in L^\infty(\Omega \times (0,T))^n$ and $A$ is a $n\times n$ symmetric positive-definite matrix with $C^2(\overline{\Omega} \times [0,T])$ coefficients, see \cite{BP}. Here, we also deal with the problem of determining an observation estimate in the general case of domain but for a linear heat equation with time-dependent coefficients
\begin{equation*}
\partial_t v - p(t)\Delta v = 0  ~~\text{in} ~~\Omega \times (0,T)
\end{equation*}
where $p\in C^1(0,T)$. 
   In this section, we will study two results of observation estimates in two different geometrical cases: The general case (Theorem \ref{nonconvex}) and the special case (Theorem \ref{convex}) when $\Omega$ is convex or star-shaped with respect to some $x_0$ such that $B(x_0,r):=\{x; |x-x_0|<r\} \subset \omega, 0<r<R:=\max\limits_{x\in \overline{\Omega}}|x-x_0|$. For the special case, we make a careful evaluation of the constants which can be explicitly computed.
First of all, we state an observation result in general case of domain $\Omega$.
\begin{theorem}
There exist constants $K=K(\Omega, \omega, p)>0$ and $\mu =\mu(\Omega, \omega, p)\in (0;1)$ such that the solution of
 \begin{equation}
\left\{\begin{array}{ll}
\partial_{t}v-p(t)\Delta v=0 & \text{in}~\Omega\times (0,T) \text{,}\\
v=0 & \text{on}~\partial \Omega \times (0,T) \text{,}\\
v(\cdot,0)\in L^2(\Omega) \text{,}
\label{J}	
\end{array}\right.
\end{equation}
satisfies
\begin{equation}
 \Vert v(\cdot,T)\Vert_{L^2(\Omega)} \leq Ke^{\frac{K}{T}} \Vert v(\cdot,T)\Vert_{L^2(\omega)}^{\mu}\Vert v(\cdot,0)\Vert_{L^2(\Omega)}^{1-\mu}\text{.}
 \label{holder}
 \end{equation}
\label{nonconvex}
\end{theorem}
\begin{corollary} For any $\varepsilon >0$, there exist positive constants $c_1$ and  $c_2$ depending on $\Omega, \omega$ and $p$ such that the following estimate holds
\begin{equation}
 \Vert v(\cdot,T)\Vert^2_{L^2(\Omega)} \leq  {c_1}e^{\frac{c_1}{T}}\frac{1}{\varepsilon^{c_2}}\Vert v(\cdot,T)\Vert^2_{L^2(\omega)}+\varepsilon\Vert v(\cdot,0)\Vert^2_{L^2(\Omega)}\text{.}
 \label{constantc}
 \end{equation}
 \label{colyoung}
\end{corollary}
\textbf{Proof of Corollary \ref{colyoung}}
\begin{proof}
It implies  from (\ref{holder}) in Theorem \ref{nonconvex} that
\begin{equation*}
 \Vert v(\cdot,T)\Vert^2_{L^2(\Omega)} \leq K^2e^{\frac{2K}{T}} \Vert v(\cdot,T)\Vert_{L^2(\omega)}^{2\mu}\Vert v(\cdot,0)\Vert_{L^2(\Omega)}^{2(1-\mu)}\text{.}
 \end{equation*}
Applying the Young's inequality $ab\leq \frac{a^m}{m}+\frac{b^q}{q}$ with
\begin{equation*} a= \left(K^{\frac{1}{\mu}}e^{\frac{K}{\mu T}}\Vert v(\cdot, T)\Vert_{L^2(\omega)}\frac{1}{\varepsilon^{\frac{1-\mu}{2\mu}}}(1-\mu)^{\frac{1-\mu}{2\mu}}\right)^{2\mu}\text{,}\end{equation*}
\begin{equation*} b= \left(\varepsilon^{\frac{1}{2}} \left(\frac{1}{1-\mu}\right)^{\frac{1}{2}}\Vert v(\cdot, 0)\Vert_{L^2(\Omega)}\right)^{2(1-\mu)}\text{,}\end{equation*}
\begin{equation*} m = \frac{1}{\mu} ~~~~~\text{and}~~~~~ q = \frac{1}{1-\mu}\text{,}\end{equation*}
we get
\begin{equation*}
\Vert v(\cdot,T)\Vert^2_{L^2(\Omega)} \leq  \mu K^{\frac{2}{\mu}}e^{\frac{2K}{\mu T}}\frac{1}{\varepsilon^{\frac{1-\mu}{\mu}}}(1-\mu)^{\frac{1-\mu}{\mu}}\Vert v(\cdot,T)\Vert^2_{L^2(\omega)}+\varepsilon\Vert v(\cdot,0)\Vert^2_{L^2(\Omega)}\text{.}
\end{equation*}
Therefore, we obtain the estimate (\ref{constantc}) with 
\begin{equation}
c_1:=\max\left\{\mu K^{\frac{2}{\mu}}(1-\mu)^{\frac{1-\mu}{\mu}}, \frac{2K}{\mu}\right\}~~\text{and} ~~c_2:=\frac{(1-\mu)}{\mu}\text{.}
\label{c1c2}
\end{equation}
\end{proof}
Our next theorem will provide us an observation result in a special geometric case with specific constants.
\begin{theorem}
Let $x_0 \in \Omega$ and $R:=\max\limits_{x\in \overline{\Omega}}|x-x_0|$. Suppose all the following assumptions hold:
\begin{enumerate}[(i)]
\item $\Omega$ is convex or star-shaped domain with respect to $x_0$,
\item $
R^2 < \frac{2p_1^2}{|p'|_\infty} ~~\text{if} ~~p \not\equiv constant
\label{small}$
 where $ |p'|_{\infty} = \sup\limits_{t\in[0,T]} |p'(t)|$;
\end{enumerate}
then the solution of (\ref{J}) satisfies (\ref{holder}) with
\[ \omega = \{x; |x-x_0|<r\}~~ \text{where} ~~0<r<R\text{,}\]
\[
 K = \max\left\{\left(4^{1+C_0(1+S_\ell)}(1+\ell)^{n+2C_0(1+S_\ell)}e^{2C_1(1+S_\ell)}e^{\frac{r^2\ell}{4p_1}}\right)^{\frac1{2(1+S_\ell)}}, \frac{r^2\ell}{4p_1(1+S_\ell)}\right\}
 \]
 and
 \[
\mu = \frac{1}{2(1+S_\ell)}\text{.}
 \]
Here \[C_0:=\frac{R^2|p'|_\infty}{2p_1^2}\text{,}\] \[C_1:=(2+n)\frac {|p'|_\infty}{p_1}\text{,}\]
\begin{equation*}
\ell :=\left\{ \begin{array}{ll}
\left(\frac{2^{2+\xi} R^2e^{C_1}}{\xi \ln \frac{3}{2}r^2}\right)^{\frac{1}{1-\xi}}-1 ~~\forall \xi \in (0,1)~~&\text{if}~~~C_0 = 0\text{,} \\
\left(\frac{4R^2e^{C_1}}{r^2\left(1-\left(\frac{2}{3}\right)^{C_0}\right)}\right)^{\frac{1}{1-C_0}}-1 ~~&\text{if}~~~C_0 > 0 
\end{array}
\right.
\end{equation*}
 and \[S_\ell := e^{C_1}\left\{
\begin{array}{ll}
\frac{\ln(1+\ell )}{\ln\frac{3}{2}}&\text{ if}~~ C_0=0\text{,}\\
\frac{(1+\ell )^{C_0}}{1-\left(\frac{2}{3}\right)^{C_0}} &\text{if}~~ C_0>0\text{.}
\end{array}\right.\]
 \label{convex}
\end{theorem}
\begin{remark} In the special case when $p \equiv 1$, the observation estimate (\ref{holder}) can be written as 
\begin{eqnarray*}
 \Vert v(\cdot,T)\Vert_{L^2(\Omega)} &\leq& \left(4(1+\ell)^{n}e^{\frac{r^2\ell}{4}\left(1+\frac 1T\right)}\Vert v(\cdot,T)\Vert_{L^2(\omega)}\right)^{\frac{1}{2\left(1+\frac{\ln(1+\ell)}{\ln \frac 32}\right)}}\Vert v(\cdot,0)\Vert_{L^2(\Omega)}^{\frac{1+2\frac{\ln(1+\ell)}{\ln \frac 32}}{2\left(1+\frac{\ln(1+\ell)}{\ln \frac 32}\right)}}
  \end{eqnarray*}
  where $\ell:=\left(\frac{2^{2+\xi} R^2}{\xi \ln \frac{3}{2}r^2}\right)^{\frac{1}{1-\xi}}-1 >1$ for any $\xi \in (0,1)$.\\
 The interested readers can compare this result with Proposition 2.1 in \cite{PW1}, Proposition 2.2 in \cite{PW2} or Theorem 4.2 in \cite{BP}.
 \end{remark}
 The main idea of the proof of both theorems is based on the logarithm convexity method (see \cite{Ve}). In order to check a kind of logarithm convexity for a suitable functional, it requires that some boundary terms must be dropped or have a good sign. This is possible under the assumption $(ii)$ in Theorem \ref{convex}. But for the general case (Theorem \ref{nonconvex}), we need a local star-shaped assumption (to get a good sign of boundary terms) and a suitable cut-off function (to drop some boundary terms). Then, thanks to the covering argument and the propagation of smallness, we get the global desired result. First of all, we need some preliminary results in the first subsection. Then, the proof of Theorem \ref{convex} and Theorem \ref{nonconvex} will be devoted in two next subsections, respectively.
\subsection{Preliminary results}
The strategy of the proof of Theorem \ref{nonconvex} and Theorem \ref{convex} consists on choosing a suitable function whose logarithm can be a convex function and considering the differential inequalities associated to this function (see Lemma \ref{lemma1}). Then by choosing a suitable weight function inspired by the heat kernel (see Corollary \ref{col1}) and solving ODE inequalities (see Lemma \ref{lemma2}), we obtain a H\"{o}lder type inequality (see Corollary \ref{col2}). The localization process in the proof of general case makes appear the function $F$ in Corollary \ref{col1}, which will be treated due to the technical Lemma \ref{lemma3}. 
\begin{lemma} 
Let $\vartheta$ be an open set in $\mathbb{R}^n$, $x_0 \in \vartheta$, $z \in H^1(0,T; H_0^1(\vartheta))$ and $\phi \in C^2(\overline{\Omega} \times (0,T))$. We define two functions from $[0,T]$ on $(0,+\infty)$ by
\begin{equation*}
y(t):=\int_{\vartheta} |z(x,t)|^2e^{\phi(x,t)}dx \text{,}
\end{equation*}
\begin{equation*}
N(t):=p(t)\frac{\int_{\vartheta} |\nabla z(x,t)|^2e^{\phi(x,t)}dx}{\int_{\vartheta} |z(x,t)| ^2e^{\phi(x,t)}dx}\text{.}
\end{equation*}
With the notations $\mathcal{G}\phi:= \partial_t \phi +p(t)\Delta \phi + p(t) |\nabla \phi|^2$ and $w:=\partial_t z - p(t)\Delta z$, the following assertions hold for any times $t>0$:
\begin{enumerate}[i/]
\item \begin{equation*}
y'(t) + 2N(t)y(t) = \int_\vartheta \mathcal{G}\phi (x,t) |z(x,t)|^2e^{\phi(x,t)}dx+2\int_\vartheta w(x,t)z(x,t)e^{\phi(x,t)}dx\text{,}
\label{lemma1i}\end{equation*}
\item \begin{eqnarray*}
N'(t)
&\leq&\frac{p'(t)}{p(t)} N (t)+ \frac{p(t)^2}{y(t)}\int_{\partial\vartheta} |\nabla z(x,t)|^2\partial_\nu\phi(x,t) e^{\phi(x,t)}dx \notag\\
&~&+\frac{p(t)}{y(t)}\int_\vartheta |\nabla z(x,t)|^2\mathcal{G}\phi(x,t)e^{\phi(x,t)}dx+\frac{1}{2y(t)}\int_\vartheta |w(x,t)|^2e^{\phi(x,t)}dx\notag\\
&~& -\frac{2p(t)^2}{y(t)} \int_{\vartheta} \nabla z(x,t) \nabla^2 \phi(x,t) \nabla z(x,t)e^{\phi(x,t)}dx\notag\\
&~&- \frac{p(t)}{y(t)^2}\int_\vartheta \mathcal{G}\phi(x,t)|z(x,t)|^2e^{\phi(x,t)}dx\int_{\vartheta}|\nabla z(x,t)|^2e^{\phi(x,t)}dx
\label{lemma1ii}
\end{eqnarray*}
\end{enumerate}
where $\nu$ is the unit outward normal vector to $\partial \vartheta$ and $\nabla^2 \phi$ is the Hessian matrix of $\phi$.
\label{lemma1}
\end{lemma}
\textbf{Proof of Lemma \ref{lemma1}}
\begin{proof} First of all, we will prove the assertion $i/$.\\
We have
\begin{equation*}
y'(t)=2\int_{\vartheta}z(x,t)\partial_t z(x,t)e^{\phi(x,t)} dx+\int_{\vartheta}|z(x,t)|^2 \partial_t\phi(x,t)e^{\phi(x,t)}dx\text{.}
\label{1}
\end{equation*}
With $w:=\partial_t z - p(t)\Delta z$, one has
\begin{eqnarray}
y'(t)&=&2\int_{\vartheta}z(x,t)w(x,t)e^{\phi(x,t)} dx+2p(t)\int_{\vartheta}z(x,t)\Delta z(x,t)e^{\phi(x,t)} dx \notag\\
&~&+\int_{\vartheta}|z(x,t)|^2 \partial_t\phi(x,t)e^{\phi(x,t)}dx\text{.}
\label{2}
\end{eqnarray}
Let us compute the second term of (\ref{2}) by using integration by parts:
\begin{eqnarray}
&~&2p(t)\int_{\vartheta}z(x,t)\Delta z(x,t)e^{\phi(x,t)} dx \notag\\
&=& -2p(t)\int_{\vartheta}  |\nabla z(x,t)|^2e^{\phi(x,t)}dx -2p(t)\int_{\vartheta} z(x,t) \nabla z(x,t) \nabla\phi(x,t) e^{\phi(x,t)}dx\notag\\
&=& -2p(t)\int_{\vartheta} |\nabla z(x,t)|^2e^{\phi(x,t)}dx-p(t)\int_{\vartheta} \nabla(|z(x,t)|^2) \nabla \phi(x,t) e^{\phi(x,t)}dx
\text{.}\label{bum}
\end{eqnarray}
We use the fact that $2z\nabla z= \nabla(|z|^2)$ to get the second equality. Integrating by parts the second term in (\ref{bum}) gives
\begin{eqnarray}
&~& -p(t)\int_{\vartheta} \nabla(|z(x,t)|^2) \nabla \phi(x,t) e^{\phi(x,t)}dx\notag\\
&=&p(t)\int_{\vartheta} |z(x,t)|^2\Delta \phi(x,t) e^{\phi(x,t)}dx + p(t)\int_{\vartheta} |z(x,t)|^2|\nabla \phi(x,t)|^2 e^{\phi(x,t)}dx \text{.}
\label{3}
\end{eqnarray}
Combining (\ref{2}) and (\ref{3}), we obtain:
\begin{eqnarray*}
y'(t) &=& -2p(t)\int_{\vartheta} |\nabla z(x,t)|^2e^{\phi(x,t)}dx+p(t)\int_{\vartheta} |z(x,t)|^2\Delta \phi(x,t) e^{\phi(x,t)}dx \notag\\
&~&+p(t)\int_{\vartheta} |z(x,t)|^2|\nabla \phi(x,t)|^2 e^{\phi(x,t)}dx+\int_{\vartheta}|z(x,t)|^2 \partial_t\phi(x,t)e^{\phi(x,t)}dx\notag\\
&~& +2\int_{\vartheta}z(x,t)w(x,t)e^{\phi(x,t)} dx\text{.}
\end{eqnarray*}
Thus, we can get the assertion i/. Now, we move to next step with the proof of assertion $ii/$.\\

\textbf{Step 1:} Compute $\frac{d}{dt}\left(p(t)\int_{\vartheta}|\nabla z(x,t)|^2e^{\phi(x,t)}dx\right)$.
\begin{eqnarray}
&~&\frac{d}{dt}\left(p(t)\int_{\vartheta}|\nabla z(x,t)|^2e^{\phi(x,t)}dx\right)\notag\\
&=& p'(t)\int_{\vartheta}|\nabla z(x,t)|^2e^{\phi(x,t)}dx+2p(t) \int_{\vartheta} \nabla z(x,t) \partial_t (\nabla z(x,t))e^{\phi(x,t)}dx \notag\\
&~& +p(t)\int_{\vartheta}|\nabla z(x,t)|^2\partial_t \phi(x,t) e^{\phi(x,t)}dx\notag\\
&=& P_1 + P_2 + P_3
\label{3*}
\end{eqnarray}
where $P_i(i=1,2,3)$ is the $i^{th}$ term in the right-hand side of (\ref{3*}). 
For the second term $P_2$, we use integration by parts, with the note that $\partial_t z = 0$ on $\partial \vartheta$, to get:
\begin{eqnarray}
P_2&=&2p(t)\int_{\vartheta}\nabla z(x,t)\nabla(\partial_t z(x,t))e^{\phi(x,t)}dx\notag\\
&=&- 2p(t)\int_{\vartheta} \Delta z(x,t) \partial_t z(x,t)e^{\phi(x,t)}dx-2p(t)\int_{\vartheta} \nabla z(x,t) \partial_t z(x,t)\nabla\phi(x,t) e^{\phi(x,t)}dx\notag\\
&=& -2 \int_{\vartheta}| \partial_t z(x,t)|^2e^{\phi(x,t)} dx+ 2\int_{\vartheta} w(x,t) \partial_t z(x,t) e^{\phi(x,t)}dx\notag\\
&~&-2p(t)  \int_{\vartheta} \partial_t z(x,t)\nabla z(x,t)\nabla \phi(x,t) e^{\phi(x,t)}dx \text{.}
\label{4}
\end{eqnarray}
The last equality is implied from the fact: $p(t)\Delta z = \partial_t z-w$. 
For the third term $P_3$, since $\mathcal{G}\phi := \partial_t \phi +p(t)\Delta \phi + p(t) |\nabla \phi|^2$, we get
\begin{eqnarray}
P_3&=&p(t)\int_{\vartheta}|\nabla z(x,t)|^2\partial_t \phi(x,t) e^{\phi(x,t)}dx\notag\\
&=&p(t)\int_{\vartheta}|\nabla z(x,t)|^2\mathcal{G}\phi(x,t)e^{\phi(x,t)}dx - p(t)^2\int_{\vartheta}|\nabla z(x,t)|^2\Delta\phi(x,t) e^{\phi(x,t)}dx\notag\\
&~& -p(t)^2\int_{\vartheta}|\nabla z(x,t)|^2|\nabla \phi(x,t)|^2 e^{\phi(x,t)}dx\text{.}
\label{5}
\end{eqnarray}
Integrating by parts the second term in (\ref{5}) gives
\begin{eqnarray}
&~&- p(t)^2\int_{\vartheta}|\nabla z(x,t)|^2\Delta\phi(x,t) e^{\phi(x,t)}dx\notag\\
&=& p(t)^2\int_{\vartheta}\nabla(|\nabla z(x,t)|^2)\nabla\phi(x,t) e^{\phi(x,t)}dx+p(t)^2\int_{\vartheta}|\nabla z(x,t)|^2|\nabla\phi(x,t)|^2 e^{\phi(x,t)}dx\notag\\
&~&- p(t)^2\int_{{\partial\vartheta}}|\nabla z(x,t)|^2\partial_\nu\phi(x,t) e^{\phi(x,t)}dx\text{.}
\label{55}
\end{eqnarray}
Now, we compute the first term in (\ref{55}) by using standard summation notations
\begin{eqnarray}
&~&p(t)^2\int_{\vartheta}\nabla(|\nabla z(x,t)|^2)\nabla\phi(x,t)e^{\phi(x,t)}dx\notag\\
&=& p(t)^2\int_{\vartheta} \partial_i (|\partial_j z(x,t)|^2)\partial_i\phi(x,t)e^{\phi(x,t)}dx\notag\\
&=& 2p(t)^2\int_{\vartheta} \partial_j z(x,t) \partial^2_{ij} z(x,t) \partial_i \phi(x,t) e^{\phi(x,t)}dx\notag\\
&=& -2p(t)^2 \int_{\vartheta} \partial^2_{jj} z(x,t) \partial_i z(x,t) \partial_i \phi(x,t)e^{\phi(x,t)}dx
-2p(t)^2 \int_{\vartheta} \partial_j z(x,t) \partial_i z(x,t)\partial^2_{ij} \phi(x,t)e^{\phi(x,t)}dx\notag\\
&~&-2p(t)^2\int_{\vartheta} \partial_j z(x,t)\partial_i z(x,t)\partial_i\phi(x,t)\partial_j\phi(x,t) e^{\phi(x,t)} dx\notag\\
&~&+2p(t)^2\int_{\partial\vartheta }\partial_j z(x,t) \partial_i z(x,t)\partial_i \phi(x,t)\nu_je^{\phi(x,t)} dx\text{.}
\end{eqnarray}
Thus, we can write
\begin{eqnarray}
&~&p(t)^2\int_{\vartheta}\nabla\left(|\nabla z(x,t)|^2\right)\nabla\phi(x,t)e^{\phi(x,t)})dx\notag\\
&=& -2p(t)^2 \int_{\vartheta} \Delta z(x,t) \nabla z(x,t) \nabla \phi(x,t)e^{\phi(x,t)}dx-2p(t)^2 \int_{\vartheta} \nabla z(x,t) \nabla^2 \phi (x,t) \nabla z(x,t)e^{\phi(x,t)}dx\notag\\
&~&-2p(t)^2 \int_{\vartheta} |\nabla z(x,t)\nabla\phi(x,t)|^2e^{\phi(x,t)}dx+2p(t)^2 \int_{\partial\vartheta} |\nabla z(x,t)|^2 \partial_\nu \phi(x,t)e^{\phi(x,t)}dx \text{.}
\label{6}
\end{eqnarray}
Combining (\ref{5}), (\ref{55}) and (\ref{6}), the third term $P_3$ in (\ref{3*}) can be computed as
\begin{eqnarray}
P_3
&=& p(t)\int_{\vartheta} |\nabla z(x,t)|^2 \mathcal{G}\phi(x,t)e^{\phi(x,t)}dx-2p(t)^2\int_\vartheta \Delta z(x,t) \nabla z(x,t) \nabla \phi(x,t)e^{\phi(x,t)}dx\notag\\
&~& -2p(t)^2\int_\vartheta \nabla z(x,t) \nabla^2 \phi (x,t) \nabla z(x,t)e^{\phi(x,t)}dx -2p(t)^2\int_\vartheta |\nabla z(x,t)\nabla \phi(x,t)|^2e^{\phi(x,t)}\notag\\
&~&+p(t)^2 \int_{\partial \vartheta}|\nabla z(x,t)|^2\partial_\nu\phi(x,t)e^{\phi(x,t)}dx\text{.}
\label{p3}
\end{eqnarray}
Thus, from above results (\ref{4}) and (\ref{p3}), (\ref{3*}) can be written
\begin{eqnarray}
&~&\frac{d}{dt}\left(p(t)\int_{\vartheta}|\nabla z(x,t)|^2e^{\phi(x,t)}dx\right)\notag\\
&=&  p'(t)\int_{\vartheta}|\nabla z(x,t)|^2e^{\phi(x,t)}dx-2 \int_{\vartheta}| \partial_t z(x,t)|^2e^{\phi(x,t)} dx+ 2\int_{\vartheta} w(x,t) \partial_t z(x,t) e^{\phi(x,t)}dx\notag\\
&~&-2p(t)  \int_{\vartheta} \partial_t z(x,t)\nabla z(x,t)\nabla \phi(x,t) e^{\phi(x,t)}dx+
p(t)\int_{\vartheta}|\nabla z(x,t)|^2\mathcal{G}\phi(x,t)e^{\phi(x,t)}dx\notag\\
&~&-2p(t)^2 \int_{\vartheta} \Delta z(x,t) \nabla z(x,t) \nabla \phi(x,t)e^{\phi(x,t)}dx-2p(t)^2 \int_{\vartheta} \nabla z(x,t) \nabla^2 \phi (x,t) \nabla z(x,t)e^{\phi(x,t)}dx\notag\\
&~&-2p(t)^2 \int_{\vartheta} |\nabla z(x,t)\nabla\phi(x,t)|^2e^{\phi(x,t)}dx+p(t)^2 \int_{\partial\vartheta} |\nabla z(x,t)|^2 \partial_\nu \phi(x,t)e^{\phi(x,t)}dx \text{.}
\label{5417}
\end{eqnarray}
Since $p(t)\Delta z = \partial_t z -w$, one has 
\begin{eqnarray*}
&~&-2p(t)^2 \int_{\vartheta} \Delta z(x,t) \nabla z(x,t) \nabla \phi(x,t)e^{\phi(x,t)}dx\notag\\
&=&-2p(t)\int_\vartheta \partial_t z(x,t)\nabla z(x,t)\nabla \phi(x,t)e^{\phi(x,t)}dx+2p(t)\int_\vartheta w(x,t)\nabla z(x,t)\nabla \phi(x,t)e^{\phi(x,t)}dx\text{.}
\end{eqnarray*}
Moreover, we also have
\begin{eqnarray}
&~&-2 \int_{\vartheta}| \partial_t z(x,t)|^2e^{\phi(x,t)} dx+ 2\int_{\vartheta} w(x,t) \partial_t z(x,t) e^{\phi(x,t)}dx\notag\\
&~&-4p(t)  \int_{\vartheta} \partial_t z(x,t)\nabla z(x,t)\nabla \phi(x,t) e^{\phi(x,t)}dx+2p(t) \int_\vartheta w(x,t)\nabla z(x,t)\nabla \phi(x,t)e^{\phi(x,t)}dx\notag\\
&~&-2p(t)^2 \int_{\vartheta} |\nabla z(x,t)\nabla\phi(x,t)|^2e^{\phi(x,t)}dx\notag\\
&=& -2\int_{\vartheta}\left(\partial_t z(x,t) +p(t)\nabla z(x,t)\nabla\phi(x,t)-\frac{1}{2}w(x,t)\right)^2e^{\phi(x,t)}dx+\frac{1}{2}\int_\vartheta |w(x,t)|^2e^{\phi(x,t)}dx\text{.}\notag\\
\label{6417}
\end{eqnarray}
Thus, (\ref{5417}) and (\ref{6417}) imply that
\begin{eqnarray}
&~&\frac{d}{dt}\left(p(t)\int_{\vartheta}|\nabla z(x,t)|^2e^{\phi(x,t)}dx\right)\notag\\
&=&  p'(t)\int_{\vartheta}|\nabla z(x,t)|^2e^{\phi(x,t)}dx+
p(t)\int_{\vartheta}|\nabla z(x,t)|^2\mathcal{G}\phi e^{\phi(x,t)}dx\notag\\
&~&-2p(t)^2 \int_{\vartheta} \nabla z(x,t) \nabla^2 \phi (x,t) \nabla z(x,t)e^{\phi(x,t)}dx+p(t)^2 \int_{\partial\vartheta} |\nabla z(x,t)|^2 \partial_\nu \phi(x,t)e^{\phi(x,t)}dx\notag\\
&~&-2 \int_{\vartheta}\left(\partial_t z(x,t) +p(t)\nabla z(x,t)\nabla\phi(x,t)-\frac{1}{2}w(x,t)\right)^2e^{\phi(x,t)}dx+\frac{1}{2}\int_\vartheta |w(x,t)|^2e^{\phi(x,t)}dx\text{.}\notag\\
\label{s1}
\end{eqnarray}

\textbf{Step 2:} Compute $y'(t)p(t)\int_{\vartheta}|\nabla z(x,t)|^2e^{\phi(x,t)}dx$.\\
From the result $i/$, we have
\begin{eqnarray}
&~&y'(t)p(t)\int_{\vartheta}|\nabla z(x,t)|^2e^{\phi(x,t)}dx\notag\\
&=& -2\left[p(t)\int_\vartheta|\nabla z(x,t)|^2e^{\phi(x,t)}dx\right]^2+2p(t)\int_\vartheta z(x,t)w(x,t)e^{\phi(x,t)}dx\int_{\vartheta}|\nabla z(x,t)|^2e^{\phi(x,t)}dx\notag\\
&~&+p(t)\int_\vartheta \mathcal{G}\phi(x,t)|z(x,t)|^2e^{\phi(x,t)}dx\int_{\vartheta}|\nabla z(x,t)|^2e^{\phi(x,t)}dx\notag\\
&=& 2A(-A+B)+p(t)\int_\vartheta \mathcal{G}\phi(x,t)|z(x,t)|^2e^{\phi(x,t)}dx\int_{\vartheta}|\nabla z(x,t)|^2e^{\phi(x,t)}dx\text{.}
\label{7417}
\end{eqnarray}
Here \[A:=p(t)\int_\vartheta|\nabla z(x,t)|^2e^{\phi(x,t)}dx\] and \[B:=\int_\vartheta z(x,t)w(x,t)e^{\phi(x,t)}dx\]
Our target is making appear the term $\partial_tz(x,t)+p(t)\nabla z(x,t) \nabla \phi(x,t)-\frac{1}{2}w(x,t)$. First of all, we compute $A$ by integrating by parts
\begin{eqnarray}
A&=&p(t)\int_{\vartheta}|\nabla z(x,t)|^2e^{\phi(x,t)}dx\notag\\
&=& -p(t)\int_\vartheta \Delta z(x,t)z(x,t)e^{\phi(x,t)}dx -p(t)\int_\vartheta \nabla z(x,t)z(x,t)\nabla \phi(x,t) e^{\phi(x,t)}dx\notag\\
&=& \int_\vartheta w(x,t)z(x,t)e^{\phi(x,t)}dx - \int_\vartheta \partial_t z(x,t)z(x,t)e^{\phi(x,t)}dx\notag\\
&~& -p(t)\int_\vartheta \nabla z(x,t)z(x,t)\nabla \phi(x,t) e^{\phi(x,t)}dx\notag\\
&=& -\int_\vartheta \left(\partial_tz(x,t)+p(t)\nabla z(x,t) \nabla \phi(x,t)-\frac{1}{2}w(x,t)\right)z(x,t)e^{\phi(x,t)}dx\notag\\
&~& +\frac{1}{2}\int_\vartheta w(x,t)z(x,t)e^{\phi(x,t)}dx\text{.}
\label{8417}
\end{eqnarray}
Thus
\begin{eqnarray}
B-A&=&\int_\vartheta \left(\partial_tz(x,t)+p(t)\nabla z(x,t) \nabla \phi(x,t)-\frac{1}{2}w(x,t)\right)z(x,t)e^{\phi(x,t)}dx\notag\\
&+&\frac{1}{2}\int_\vartheta w(x,t)z(x,t)e^{\phi(x,t)}dx\text{.}
\label{9417}
\end{eqnarray}
Combining (\ref{7417}), (\ref{8417}) and (\ref{9417}), one gets
\begin{eqnarray}
&~&y'(t)p(t)\int_{\vartheta}|\nabla z(x,t)|^2e^{\phi(x,t)}dx\notag\\
&=&  \frac{1}{2}\left(\int_\vartheta w(x,t)z(x,t)e^{\phi(x,t)}dx\right)^2\notag\\
&~& - 2\left(\int_\vartheta \left(\partial_tz(x,t)+p(t)\nabla z(x,t) \nabla \phi(x,t)-\frac{1}{2}w(x,t)\right)z(x,t)e^{\phi(x,t)}dx\right)^2\notag\\
&~&+p(t)\int_\vartheta \mathcal{G}\phi(x,t)|z(x,t)|^2e^{\phi(x,t)}dx\int_{\vartheta}|\nabla z(x,t)|^2e^{\phi(x,t)}dx\text{.}\notag\\
\label{s2}
\end{eqnarray}

\textbf{Step 3:} Compute $N'(t)$.\\\\
We have
\begin{equation*}
N'(t)= \frac{1}{y(t)^2}\left(y(t)\frac{d}{dt}\left(p(t)\int_\vartheta |\nabla z(x,t)|^2e^{\phi(x,t)}\right)-y'(t)p(t)\int_\vartheta |\nabla z(x,t)|^2e^{\phi(x,t)}\right)\text{.}
\end{equation*}
The result (\ref{s1}) in Step 1 and (\ref{s2}) in Step 2 provide us
\begin{eqnarray*}
N'(t)
&=& \frac{p'(t)}{p(t)} N (t)+ \frac{p(t)^2}{y(t)}\int_{\partial\vartheta} |\nabla z(x,t)|^2\partial_\nu\phi(x,t) e^{\phi(x,t)}dx +\frac{p(t)}{y(t)}\int_\vartheta |\nabla z(x,t)|^2\mathcal{G}\phi e^{\phi(x,t)}dx\notag\\
&~& -\frac{2p(t)^2}{y(t)} \int_{\vartheta} \nabla z(x,t) \nabla^2 \phi (x,t) \nabla z(x,t)e^{\phi(x,t)}dx+\frac{1}{2y(t)}\int_\vartheta |w(x,t)|^2e^{\phi(x,t)}dx\notag\\
&~&-\frac{2}{y(t)} \int_{\vartheta}\left(\partial_t z(x,t) +p(t)\nabla z(x,t)\nabla\phi(x,t)-\frac{1}{2}w(x,t)\right)^2e^{\phi(x,t)}dx\notag\\
&~& +\frac{2}{y(t)^2}\left(\int_\vartheta \left(\partial_tz(x,t)+p(t)\nabla z(x,t) \nabla \phi(x,t)-\frac{1}{2}w(x,t)\right)z(x,t)e^{\phi(x,t)}dx\right)^2\notag\\
&~&-\frac{1}{2y(t)^2}\left(\int_\vartheta w(x,t)z(x,t)e^{\phi(x,t)}dx\right)^2\notag\\
&~& - \frac{p(t)}{y(t)^2}\int_\vartheta \mathcal{G}\phi|z(x,t)|^2e^{\phi(x,t)}dx\int_{\vartheta}|\nabla z(x,t)|^2e^{\phi(x,t)}dx\text{.}
\end{eqnarray*}
Thanks to Cauchy-Schwarz inequality:
\begin{eqnarray*}
&~&\left(\int_\vartheta \left(\partial_tz(x,t)+p(t)\nabla z(x,t) \nabla \phi(x,t)-\frac{1}{2}w(x,t)\right)z(x,t)e^{\phi(x,t)}dx\right)^2\notag\\
&\leq& \int_{\vartheta}\left(\partial_t z(x,t) +p(t)\nabla z(x,t)\nabla\phi(x,t)-\frac{1}{2}w(x,t)\right)^2e^{\phi(x,t)}dx \int_\vartheta |z(x,t)|^2e^{\phi(x,t)}dx\text{,}\notag\\
\end{eqnarray*}
we receive the assertion ii/.
\end{proof}
Now, by choosing an explicit weight function $e^\phi$ inspired from the heat kernel, we get the following result.
\begin{corollary}
Under the same assumption in Lemma \ref{lemma1}, put $R:=\max\limits_{x\in \overline{\vartheta}}|x-x_0|$. Assume that $\vartheta$ be a convex domain or star-shaped with respect to $x_0$. For any $\rho >0$, with $\phi$ is chosen as below
\begin{equation}
\phi(x,t):=\frac{-|x-x_0|^2}{4p(T)(T-t+\rho)}-\frac{n}{2}\ln(T-t+\rho)\text{,} 
\end{equation}
we obtain two following estimates:
\begin{enumerate}[i/]
\item \begin{equation}
\left|y'(t)+2N(t)y(t)\right| \leq \left(\frac{C_0}{T-t+\rho}+C_1 \right)y(t) +2\int_\vartheta |w(x,t)z(x,t)|e^{\phi(x,t)}dx\text{,}
\label{col1i}
\end{equation}
\item \begin{equation}
N'(t)\leq\left(\frac{1+C_0}{T-t+\rho}+C_1\right)N(t)+\frac{1}{2}\frac{ \int_{\vartheta} |w(x,t)|^2e^{\phi(x,t)}dx}{y(t)}  
\label{col1ii}
\end{equation}
\end{enumerate}
where \begin{equation*}C_0=\frac{|p'|_\infty R^2 }{2p_1^2}~~\text{and}~~ C_1=\left(2+n\right)\frac{|p'|_\infty}{p_1}\text{.}
\label{c0c1}
\end{equation*}
\label{col1}
\end{corollary}
\textbf{Proof of Corollary \ref{col1}}
\begin{proof}
Obviously, we can easily check the following properties of the function $\phi$:
\begin{enumerate}[(1)]
\item $\partial_t \phi+p(T)\Delta \phi + p(T) |\nabla\phi|^2 =0$,
\item $\nabla \phi = \frac{-(x-x_0)}{2p(T)(T-t+\rho)}$,
\item $\Delta \phi = \frac{-n}{2p(T)(T-t+\rho)}$,
\item $\nabla^2\phi =\frac{-1}{2p(T)(T-t+\rho)} I_n$
where $I_n$ is the identity matrix of size $n$.
\end{enumerate}
Remind that $\mathcal{G}\phi=\partial_t \phi+p(t)\Delta \phi + p(t)|\nabla \phi|^2$. Thanks to properties (1), (2) and (3), we get
\begin{eqnarray}
|\mathcal{G}\phi|&\leq&\left|p(t)-p(T)\right|\Delta\phi + \left|p(t)-p(T)\right||\nabla \phi|^2\notag\\
&\leq& \frac{n|p'|_\infty}{2p(T)}+\frac{|p'|_\infty R^2}{4p(T)^2}\frac{1}{T-t+\rho}\notag\\
&\leq& \frac{n|p'|_\infty}{2p_1}+\frac{|p'|_\infty R^2}{4p_1^2}\frac{1}{T-t+\rho}\text{.}
\label{gbig}
\end{eqnarray}
Hence, from result $i/$ in Lemma \ref{lemma1}, we get the assertion $i/$. Now, we turn to prove the assertion $ii/$. Thanks to the assumption that $\vartheta$ is star-shaped with respect to $x_0$, one has
\begin{equation}
\partial_\nu \phi = -\frac{(x-x_0)\nu}{2p(T)(T-t+\rho)} \leq 0 ~~\forall x\in  \partial \vartheta\text{.}
\label{nu}
\end{equation}
Furthermore, property (4) implies
\begin{equation}
\int_{\Omega} \nabla z(x,t) \nabla^2 \phi (x,t) \nabla z(x,t)e^{\phi(x,t)}dx = \frac{-1}{2p(T)(T-t+\rho)}\int_\Omega |\nabla z(x,t)|^2e^{\phi(x,t)}dx\text{.}
\label{ij}
\end{equation}
Consequently, combining result $ii/$ in Lemma \ref{lemma1} with (\ref{gbig}), (\ref{nu}) and (\ref{ij}), we get the assertion $ii/$.\end{proof}
Now, the following lemma will solve the ODE inequalities getting from Corollary \ref{col1}.
\begin{lemma}
Let $\rho>0$, $F\in C^0([0,T])$. Suppose two positive functions $y, N \in C^1([0,T])$ satisfy the following conditions
\begin{enumerate}
\item \begin{equation}
\left| y'(t) + 2N(t)y(t) \right| \leq \left(\frac{C_0}{T-t+\rho}+C_1+F(t)\right)y(t)  \text{,}\label{ass1}
\end{equation}
\item \begin{equation}
N'(t)\leq\left(\frac{1+C_0}{T-t+\rho}+C_1\right)N(t)+\frac 12 F(t) \label{ass2}
\end{equation}
\end{enumerate}
where $C_0, C_1 >0$. Then for any $0\leq t_1<t_2<t_3\leq T$, one has
\begin{equation}
\left(y(t_2)\right)^{1+M} \leq e^G \left(\frac{T-t_1+\rho}{T-t_3+\rho}\right)^{C_0(1+M)}y(t_3) \left(y(t_1)\right)^M 
\label{ine}
\end{equation}
with
\[
M=\frac{\int_{t_2}^{t_3}\frac{e^{C_1s}}{(T-s+h)^{1+C_0}}ds}{\int_{t_1}^{t_2}\frac{e^{C_1s}}{(T-s+h)^{1+C_0}}ds}
\text{,}\]
\[ G = (1+M)\left[(t_3-t_1)\int_{t_1}^{t_3}F(s)ds+\int_{t_1}^{t_3}F(s)ds+(t_3-t_1)C_1\right] \text{.}\]
\label{lemma2}
\end{lemma}
\textbf{Proof of Lemma \ref{lemma2}}
\begin{proof}
From (\ref{ass2}), we get:
\begin{equation}
 \left(N(t)(T-t+\rho )^{1+C_0}e^{-C_1t}\right)'\leq \frac{1}{2}F(t)(T-t+\rho )^{1+C_0}e^{-C_1t} \text{.}
 \label{N}
 \end{equation}

\underline{For $t_1<t<t_2$:}\\\\
Integrating ($\ref{N}$) over $(t;t_2)$ gives
\begin{eqnarray}
N(t) &\geq& N(t_2) \left(\frac{T-t_2+\rho }{T-t+\rho }\right)^{1+C_0}e^{-C_1(t_2-t)} \notag\\ 
&~&-\frac{1}{2}e^{C_1t}\left(\frac{1}{T-t+\rho }\right)^{1+C_0}\int_t^{t_2}F(s)(T-s+\rho )^{1+C_0}e^{-C_1s}ds\text{.}
\end{eqnarray}
Using the fact that $(T-s+\rho )^{1+C_0}e^{-C_1s}\leq (T-t+\rho )^{1+C_0}e^{-C_1t}~~\forall s\geq t$, one gets
\begin{equation}
N(t) \geq Q(t_2)\frac{e^{C_1t}}{(T-t+\rho )^{1+C_0}}-\frac{1}{2}\int_{t_1}^{t_2}F(s)ds
\label{Na}
\end{equation}
where $Q(t_2) = e^{-C_1t_2}(T-t_2+\rho )^{1+C_0}N(t_2)$. 
From (\ref{ass1}), we also have
\begin{equation}
y'(t)+2N(t)y(t) \leq \left(\frac{C_0}{T-t+\rho }+C_1+F(t)\right)y(t) \text{.}
\end{equation}
Combining to ($\ref{Na}$), we obtain:
\begin{equation*}
 y'(t) + \left(2Q(t_2)\frac{e^{C_1t}}{(T-t+\rho )^{1+C_0}}-\int_{t_1}^{t_2}F(s)ds-\frac{C_0}{T-t+\rho }-C_1-F(t)\right)y(t) \leq 0 \text{.}
\end{equation*}
It is equivalent to 
\begin{equation}
\left(y(t)e^{2Q(t_2)\int_0^t \frac {e^{C_1s}}{(T-s+\rho )^{1+C_0}}ds}e^{-\left(\int_{t_1}^{t_2}F(s)ds + C_1\right)t}(T-t+\rho )^{C_0}e^{-\int_0^tF(s)ds}\right)'\leq 0 \text{.}
\label{de}
\end{equation}
Integrating (\ref{de}) over $(t_1;t_2)$, one has
\begin{equation}
y(t_1) \geq y(t_2) e^{2Q(t_2)\int_{t_1}^{t_2} \frac {e^{C_1s}}{(T-s+\rho )^{1+C_0}}ds}e^{-\left(\int_{t_1}^{t_2}F(s)ds + C_1\right)(t_2-t_1)}\left(\frac{T-t_2+\rho }{T-t_1+\rho }\right)^{C_0}e^{-\int_{t_1}^{t_2}F(s)ds} \text{.}
\label{t1}
\end{equation}

\underline{For $t_2<t<t_3$:}\\\\
Integrating ($\ref{N}$) over $(t_2;t)$ gives
\begin{eqnarray}
N(t) &\leq& N(t_2) \left(\frac{T-t_2+\rho }{T-t+\rho }\right)^{1+C_0}e^{-C_1(t_2-t)} \notag\\ 
&~&+\frac{1}{2}e^{C_1t}\left(\frac{1}{T-t+\rho }\right)^{1+C_0}\int_{t_2}^{t}F(s)(T-s+\rho )^{1+C_0}e^{-C_1s}ds\text{.}
\end{eqnarray}
Using the fact that $(T-s+\rho )^{1+C_0}e^{-C_1s}\leq (T-t_2+\rho )^{1+C_0}e^{-C_1t_2}~~\forall s\geq t_2$, we obtain
\begin{eqnarray}
N(t) \leq Q(t_2) \left(\frac{1}{T-t+\rho }\right)^{1+C_0}e^{C_1t}+ \frac{1}{2}e^{C_1(t-t_2)}\left(\frac{T-t_2+\rho }{T-t+\rho }\right)^{1+C_0}\int_{t_2}^{t_3}F(s)ds\text{.}
\label{Na2}
\end{eqnarray}
From (\ref{ass1}), we also have
\begin{equation}
y'(t)+2N(t)y(t) \geq -\left(\frac{C_0}{T-t+\rho }+C_1+F(t)\right)y(t) \text{.}
\label{sunday}
\end{equation}
It deduces from ($\ref{Na2}$) and (\ref{sunday}) that
\begin{eqnarray}
y'(t)+\left(2Q(t_2)\frac{e^{C_1t}}{(T-t+\rho )^{1+C_0}}+e^{C_1(t-t_2)}\left(\frac{T-t_2+\rho }{T-t+\rho }\right)^{1+C_0}\int_{t_2}^{t_3}F(s)ds\right)\notag\\
\geq-\left(\frac{C_0}{T-t+\rho }+C_1+F(t)\right)y(t) \text{.}
\end{eqnarray}
It is equivalent to 
\begin{equation}
\left(y(t)e^{2Q(t_2)\int_0^t \frac {e^{C_1s}}{(T-s+\rho )^{1+C_0}}ds+\int_{t_2}^{t_3}F(s)ds\int_0^te^{C_1(s-t_2)}\left(\frac{T-t_2+\rho }{T-s+\rho }\right)^{1+C_0}ds+C_1t+\int_0^tF(s)ds}\frac{1}{(T-t+\rho )^{C_0}}\right)'\geq 0 \text{.}
\label{de2}
\end{equation}
Integrating ($\ref{de2}$) over $(t_2;t_3)$ gives
\begin{eqnarray}
y(t_2) \leq y(t_3) e^{2Q(t_2)\int_{t_2}^{t_3} \frac {e^{C_1s}}{(T-s+\rho )^{1+C_0}}ds}e^{\int_{t_2}^{t_3}F(s)ds\int_{t_2}^{t_3}e^{C_1(s-t_2)}\left(\frac{T-t_2+\rho }{T-s+\rho }\right)^{1+C_0}ds}\notag\\
\times \left(\frac{T-t_2+\rho }{T-t_3+\rho }\right)^{C_0}e^{C_1(t_3-t_2)}e^{\int_{t_2}^{t_3}F(s)ds} \text{.}
\end{eqnarray}
Combining to ($\ref{t1}$), one gets
\begin{eqnarray}
y(t_2) &\leq& y(t_3) \left(\frac{y(t_1)}{y(t_2)}\right)^{M}\left(\frac {T-t_1+\rho }{T-t_2+\rho }\right)^{MC_0}e^{\left(\int_{t_1}^{t_2}F(s)ds+C_1\right)M(t_2-t_1)}e^{M\int_{t_1}^{t_2}F(s)ds}\notag\\
&~& \times e^{C_1(t_3-t_2)}e^{\int_{t_2}^{t_3}F(s)ds\int_{t_2}^{t_3}e^{C_1(s-t_2)}\left(\frac{T-t_2+\rho }{T-s+\rho }\right)^{1+C_0}ds}e^{\int_{t_2}^{t_3}F(s)ds}\left(\frac{T-t_2+\rho }{T-t_3+\rho }\right)^{C_0}\notag\\
\label{bumbum}
\end{eqnarray}
where
\begin{equation}
M=\frac{\int_{t_2}^{t_3}\frac{e^{C_1s}}{(T-s+\rho )^{1+C_0}}ds}{\int_{t_1}^{t_2}\frac{e^{C_1s}}{(T-s+\rho )^{1+C_0}}ds} \text{.}
\end{equation}
We also have
\begin{eqnarray}
\int_{t_2}^{t_3}e^{C_1(s-t_2)}\left(\frac{T-t_2+\rho }{T-s+\rho }\right)^{1+C_0}ds&=& \frac{\int_{t_2}^{t_3}\frac{e^{C_1s}}{(T-s+\rho )^{1+C_0}}ds}{\frac{e^{C_1t_2}}{(T-t_2+\rho )^{1+C_0}}}\notag\\
&\leq& \left(t_2 - t_1\right)\frac{\int_{t_2}^{t_3}\frac{e^{C_1s}}{(T-s+\rho )^{1+C_0}}ds}{\int_{t_1}^{t_2}\frac{e^{C_1s}}{(T-s+\rho )^{1+C_0}}ds}\notag\\
&\leq& (t_2-t_1)M \text{.}
\end{eqnarray}
Hence, it is deduced from (\ref{bumbum}) that
\begin{eqnarray}
&~&\left(y(t_2)\right)^{1+M}\notag\\
&\leq& y(t_3)\left(y(t_1)\right)^Me^{M(t_2-t_1)\int_{t_1}^{t_2}F(s)ds}e^{C_1M(t_2-t_1)}e^{M\int_{t_1}^{t_2}F(s)ds}\notag\\
&~& \times e^{M(t_2-t_1)\int_{t_2}^{t_3}F(s)ds}e^{C_1(t_3-t_2)}e^{\int_{t_2}^{t_3}F(s)ds}\notag\\
&~& \times \left(\frac{T-t_2+\rho }{T-t_3+\rho }\right)^{C_0}\left(\frac{T-t_1+\rho }{T-t_2+\rho }\right)^{MC_0}\notag\\
&\leq&y(t_3)\left(y(t_1)\right)^Me^{(1+M)\left[(t_3-t_1)\int_{t_1}^{t_3}F(s)ds+\int_{t_1}^{t_3}F(s)ds+(t_3-t_1)C_1\right]}\left(\frac{T-t_1+\rho }{T-t_3+\rho }\right)^{(1+M)C_0}\text{.}\notag\\
\end{eqnarray}
\end{proof}
We move to an application of this Lemma with specific choice of time.
\begin{corollary}
Under the assumption of Lemma \ref{lemma2}, for any $\rho>0$ and $\ell>1$ such that $\ell \rho \leq \min\{\frac{1}{2};\frac{T}{4}\}$, one has
\begin{eqnarray}
\left(y(T-\ell \rho)\right)^{1+M_\ell } \leq e^{G_\ell } (1+2\ell )^{2C_0(1+M_\ell )} y(T)\left(y(T-2\ell \rho)\right)^{M_\ell }
\end{eqnarray}
where 
\begin{equation}
M_\ell=\frac{\int_{T-\ell\rho }^{T}\frac{e^{C_1s}}{(T-s+\rho )^{1+C_0}}ds}{\int_{T-2\ell\rho }^{T-\ell\rho }\frac{e^{C_1s}}{(T-s+\rho )^{1+C_0}}ds} \text{.}
\end{equation}
and
\begin{equation}G_\ell  = (1+M_\ell )\left( 2\int_{T-2\ell \rho}^{T}F(s)ds + C_1\right)\text{.}\end{equation}
Moreover, the upper bound of $M_\ell$ can be given as
\begin{equation*}M_\ell \leq S_\ell := e^{C_1} \left\{
\begin{array}{ll}
\frac{\ln(1+\ell )}{\ln\frac{3}{2}}&\text{ if}~~ C_0=0\\
\frac{(1+\ell )^{C_0}}{1-\left(\frac{2}{3}\right)^{C_0}} &\text{if}~~ C_0>0
\end{array} \right.\text{.} \end{equation*}
\label{col2}
\end{corollary}
\textbf{Proof of Corollary \ref{col2}}
\begin{proof}
Now, for $\ell >1$ and $\ell\rho <\min\{\frac{1}{2}; \frac{T}{4}\}$, applying Lemma \ref{lemma2} for $t_1= T-2\ell\rho ; ~ t_2= T-\ell\rho ; ~ t_3=T$ , we get
\begin{eqnarray}
\left(y(T-\ell\rho )\right)^{1+M_\ell}\leq y(T)\left(y(T-2\ell\rho )\right)^{M_\ell}e^{(1+M_\ell)\left(2\int^T_{T-2\ell\rho }F(s)ds+C_1\right)}(1+2\ell)^{(1+M_\ell)C_0}
\label{step4}
\end{eqnarray}
with
\begin{equation*}
M_\ell=\frac{\int_{T-\ell\rho }^{T}\frac{e^{C_1s}}{(T-s+\rho )^{1+C_0}}ds}{\int_{T-2\ell\rho }^{T-\ell\rho }\frac{e^{C_1s}}{(T-s+\rho )^{1+C_0}}ds} \text{.}
\end{equation*}
If $C_0=0$ then
\begin{equation*}
M_\ell=\frac{\int_{T-\ell \rho}^{T}\frac{e^{C_1s}}{(T-s+\rho)}ds}{\int_{T-2\ell \rho}^{T-\ell \rho}\frac{e^{C_1s}}{(T-s+\rho)}ds}
\leq e^{C_1}\frac{\ln(1+\ell)}{\ln\frac{1+2\ell}{1+\ell}}
\leq e^{C_1}\frac{\ln(1+\ell)}{\ln\frac{3}{2}}\text{.}
\end{equation*}
If $C_0> 0$ then
\begin{equation*}
M_\ell \leq e^{2\ell \rho C_1}\frac{(1+\ell)^{C_0}-1}{\frac{(1+2\ell)^{C_0}-(1+\ell)^{C_0}}{(1+2\ell)^{C_0}}}\leq e^{C_1}\frac{(1+\ell)^{C_0}}{1-\left(\frac{2}{3}\right)^{C_0}}\text{.}
\end{equation*}
\end{proof}
\begin{lemma}
Let $x_0 \in \Omega$, $\varrho > 0$ and $0<\epsilon < \frac{\varrho}{2}$. Let $v$ be the solution of (\ref{J}). 
Then there exist constants $E_1>1$, $E_2>0$ and $E_3>0$, which all depend on $\varrho$ and $\epsilon$, such that the following estimate holds
\begin{equation*} \frac{\int_\Omega \left| v(x,0)\right|^2dx}{\int_{\Omega \cap B(x_0,\varrho)} \left|v(x,t)\right|^2dx} \leq E_1e^{\frac{E_1}{\theta}} ~~\forall~~  \frac T2<T-E_2\theta\leq t \leq T\text{.}
\end{equation*}
Here
 \begin{equation*} \frac 1\theta =\ln\left(E_3e^{\frac {E_3}T}\frac{\int_\Omega \left|v(x,0)\right|^2dx}{\int_{\Omega \cap B(x_0,\varrho - 2\epsilon)}\left| v(x,T)\right|^2dx}\right)\text{.}\end{equation*}
\label{lemma3}
\end{lemma}
\textbf{Proof of Lemma \ref{lemma3}}\\
In order to get a local estimate, we need to use a weight function $e^{-\frac{|x-x_0|^2}{\hbar}}$ ($\hbar$ will be chosen later) and a cut-off function $\Psi$ on $B(x_0,\varrho)$. After finding an ODE inequality for 
$\int_{\Omega \cap B(x_0,\varrho)} \left|\Psi(x) v(x,t)\right|^2 e^{\frac{-\left|x-x_0\right|^2}{\hbar}} dx$ (see Step 1) and solving it (see Step 2), we can get the final result with a suitable choice of $\hbar$ (see Step 3). Now, we start the proof with the definition of the following cut-off function.
 \begin{proof}
Let $\Psi \in C^\infty_0(B(x_0,\varrho))$ such that 
\begin{equation*}
\left\{\begin{array}{ll}
\Psi = 1 &\text{in} ~~B(x_0,\varrho - \epsilon) \text{,}\\
 \Psi\in(0;1) &\text{in}~~B(x_0,\varrho)
 \text{.}\end{array}\right.
 \end{equation*}
 \textbf{Step 1:} For $\hbar <1$, find an ODE inequality for  $\int_{\Omega \cap B(x_0,\varrho)} \left|\Psi(x) v(x,t)\right|^2 e^{\frac{-\left|x-x_0\right|^2}{\hbar}} dx$.\\
We have
\begin{eqnarray}
&~& \frac{d}{dt}\int_{\Omega \cap B(x_0,\varrho)} \left|\Psi(x) v(x,t)\right|^2 e^{\frac{-\left|x-x_0\right|^2}{\hbar}} dx\notag\\
&=& 2\int_{\Omega \cap B(x_0,\varrho)} \left|\Psi(x)\right|^2 v(x,t) \partial_t v(x,t) e^{\frac{-\left|x-x_0\right|^2}{\hbar}} dx\notag\\
&=& 2p(t)\int_{\Omega \cap B(x_0,\varrho)} \left|\Psi(x)\right|^2 v(x,t) \Delta v(x,t) e^{\frac{-\left|x-x_0\right|^2}{\hbar}} dx
\text{.}
\label{smile1}
\end{eqnarray}
By integrating by parts with the fact that $\Psi v = 0$ on $\partial(\Omega \cap B(x_0, \varrho))$ one obtains
\begin{eqnarray}
&~& \int_{\Omega \cap B(x_0,\varrho)} \left|\Psi(x)\right|^2 v(x,t) \Delta v(x,t) e^{\frac{-\left|x-x_0\right|^2}{\hbar}} dx\notag\\
&=& - \int_{\Omega \cap B(x_0,\varrho)} \nabla v(x,t) \nabla(|\Psi(x)|^2 v(x,t)e^{\frac{-\left|x-x_0\right|^2}{\hbar}}) dx\notag\\
&=& -2 \int_{\Omega \cap B(x_0,\varrho)} \Psi(x) \nabla \Psi(x) \nabla v(x,t)v(x,t)e^{\frac{-\left|x-x_0\right|^2}{\hbar}}dx\notag\\
&~& -\int_{\Omega \cap B(x_0,\varrho)} |\Psi(x)|^2 \left|\nabla v(x,t)\right|^2e^{\frac{-\left|x-x_0\right|^2}{\hbar}}dx\notag\\
&~& + 2 \int_{\Omega \cap B(x_0,\varrho)} |\Psi(x)|^2 v(x,t) \nabla v(x,t)(x-x_0)\frac{1}{\hbar} e^{\frac{-\left|x-x_0\right|^2}{\hbar}}dx\notag\\
&:=& P_1 +P_2+P_3
\label{smile2}
\end{eqnarray}
where $P_i(i=1,2,3)$ is the $i^{th}$ term in the right-hand side of (\ref{smile2}). Now, thanks to Cauchy-Schwarz inequality, we get
\begin{eqnarray*}
P_1 &=&-2\int_{\Omega \cap B(x_0,\varrho)} \Psi(x) \nabla \Psi(x) \nabla v(x,t)v(x,t)e^{\frac{-\left|x-x_0\right|^2}{\hbar}}dx\notag\\
&\leq& 2\left(\int_{\Omega \cap B(x_0,\varrho)}\left| \Psi(x) \nabla v(x,t)\right|^2e^{\frac{-\left|x-x_0\right|^2}{\hbar}}dx\right)^{\frac{1}{2}}\left(\int_{\Omega \cap B(x_0,\varrho)}\left| v(x,t) \nabla \Psi(x)\right|^2e^{\frac{-\left|x-x_0\right|^2}{\hbar}}dx\right)^{\frac{1}{2}}\notag\\
\end{eqnarray*}
and
\begin{eqnarray*}
P_3&=&2 \int_{\Omega \cap B(x_0,\varrho)} |\Psi(x)|^2 v(x,t) \nabla v(x,t)(x-x_0)\frac{1}{\hbar} e^{\frac{-\left|x-x_0\right|^2}{\hbar}}dx\notag\\
&\leq& 2\left(\int_{\Omega \cap B(x_0,\varrho)}\left| \Psi(x) \nabla v(x,t)\right|^2e^{\frac{-\left|x-x_0\right|^2}{\hbar}}dx\right)^{\frac{1}{2}}\notag\\
&~&\times \left(\int_{\Omega \cap B(x_0,\varrho)}\left| \Psi(x) v(x,t)\frac{(x-x_0)}{\hbar}\right|^2e^{\frac{-\left|x-x_0\right|^2}{\hbar}}dx\right)^{\frac{1}{2}}
\text{.}\notag\\
\end{eqnarray*}
Thus
\begin{eqnarray}
P_1+P_3
&\leq& 2\left(\int_{\Omega \cap B(x_0,\varrho)}\left| \Psi(x) \nabla v(x,t)\right|^2e^{\frac{-\left|x-x_0\right|^2}{\hbar}}dx\right)^{\frac{1}{2}}\notag\\
& \times& \left[\left(\int_{\Omega \cap B(x_0,\varrho)}\left| v(x,t) \nabla \Psi(x)\right|^2e^{\frac{-\left|x-x_0\right|^2}{\hbar}}dx\right)^{\frac{1}{2}}\right.\notag\notag\\
&~&~~~\left.+\left(\int_{\Omega \cap B(x_0,\varrho)}\left| \Psi(x) v(x,t)\frac{(x-x_0)}{\hbar}\right|^2e^{\frac{-\left|x-x_0\right|^2}{\hbar}}dx\right)^{\frac{1}{2}}\right]\text{.}
\end{eqnarray}
It follows from $2ab \leq a^2 + b^2$ and $(a+b)^2 \leq 2a^2 + 2b^2~~\forall a,b>0$ that
\begin{eqnarray}
P_1+P_3
&\leq& \int_{\Omega \cap B(x_0,\varrho)}\left| \Psi(x) \nabla v(x,t)\right|^2e^{\frac{-\left|x-x_0\right|^2}{\hbar}}dx\notag\\
&+& \left[\left(\int_{\Omega \cap B(x_0,\varrho)}\left| v(x,t) \nabla \Psi(x)\right|^2e^{\frac{-\left|x-x_0\right|^2}{\hbar}}dx\right)^{\frac{1}{2}}\right.\notag\\
&~&~~\left.+\left(\int_{\Omega \cap B(x_0,\varrho)}\left| \Psi(x) v(x,t)\frac{(x-x_0)}{\hbar}\right|^2e^{\frac{-\left|x-x_0\right|^2}{\hbar}}dx\right)^{\frac{1}{2}}\right]^2\notag\\
&\leq& \int_{\Omega \cap B(x_0,\varrho)}\left| \Psi(x) \nabla v(x,t)\right|^2e^{\frac{-\left|x-x_0\right|^2}{\hbar}}dx + 2\int_{\Omega \cap B(x_0,\varrho)}\left| v(x,t) \nabla \Psi(x)\right|^2e^{\frac{-\left|x-x_0\right|^2}{\hbar}}dx\notag\\
&~&+2\int_{\Omega \cap B(x_0,\varrho)}\left| \Psi(x) v(x,t)\frac{(x-x_0)}{\hbar}\right|^2e^{\frac{-\left|x-x_0\right|^2}{\hbar}}dx
\label{smile3}
\end{eqnarray}
Combining (\ref{smile1}), (\ref{smile2}) and (\ref{smile3}), we can conclude
\begin{eqnarray}
&~& \frac{d}{dt}\int_{\Omega \cap B(x_0,\varrho)} \left|\Psi(x) v(x,t)\right|^2 e^{\frac{-\left|x-x_0\right|^2}{\hbar}} dx\notag\\
 &\leq& 4p(t)\int_{\Omega \cap B(x_0,\varrho)}\left| v(x,t) \nabla \Psi(x)\right|^2e^{\frac{-\left|x-x_0\right|^2}{\hbar}}dx\notag\\
 &~&+4p(t)\int_{\Omega \cap B(x_0,\varrho)}\left| \Psi(x) v(x,t)\frac{(x-x_0)}{\hbar}\right|^2e^{\frac{-\left|x-x_0\right|^2}{\hbar}}dx\notag\\
 &\leq& 4p_2\int_{\Omega \cap B(x_0,\varrho)}\left| v(x,t)\nabla \Psi(x) \right|^2e^{\frac{-\left|x-x_0\right|^2}{\hbar}}dx\notag\\
&~& +4p_2\frac{\varrho^2}{\hbar^2} \int_{\Omega \cap B(x_0,\varrho)}\left| \Psi(x) v(x,t)\right|^2e^{\frac{-\left|x-x_0\right|^2}{\hbar}}dx\text{.}
\end{eqnarray}
Moreover, due to $\nabla \Psi(\cdot) = 0$ in $\Omega \cap B(x_0, \varrho - \epsilon)$, one has
\begin{eqnarray}
\int_{\Omega \cap B(x_0,\varrho)} \left| v(x,t)\nabla \Psi(x)\right|^2 e^{\frac{-\left|x-x_0\right|^2}{\hbar}}dx 
&=& \int_{\Omega \cap \{\left|x-x_0\right|\geq \varrho-\epsilon\}} \left| v(x,t)\nabla \Psi(x)\right|^2 e^{\frac{-\left|x-x_0\right|^2}{\hbar}}dx\notag\\
&\leq&|\nabla \Psi|_\infty e^{\frac{-(\varrho-\epsilon)^2}{\hbar}}\int_\Omega \left|v(x,t)\right|^2dx\notag\\
&\leq&|\nabla \Psi|_\infty e^{\frac{-(\varrho-\epsilon)^2}{\hbar}}\int_\Omega \left|v(x,0)\right|^2dx
\text{.}\end{eqnarray}
Thus
\begin{eqnarray}
&~&\frac{d}{dt}\int_{\Omega \cap B(x_0,\varrho)} \left|\Psi(x) v(x,t)\right|^2 e^{\frac{-\left|x-x_0\right|^2}{\hbar}}dx\notag\\
&\leq& \frac{4p_2\varrho^2}{\hbar^2}\int_{\Omega \cap B(x_0,\varrho)} \left|\Psi(x) v(x,t)\right|^2 e^{\frac{-\left|x-x_0\right|^2}{\hbar}}dx+4p_2| \nabla \Psi|_\infty^2e^{\frac{-(\varrho-\epsilon)^2}{\hbar}}\int_\Omega \left|v(x,0)\right|^2dx\text{.}\notag\\
\label{ode}
\end{eqnarray}

\textbf{Step 2:} Solve ODE inequality.\\
It deduces from (\ref{ode}) that
\begin{eqnarray}
&~&\frac{d}{dt}\left(e^{-4p_2\frac{\varrho^2}{\hbar^2}t}\int_{\Omega \cap B(x_0,\varrho)} \left|\Psi(x) v(x,t)\right|^2 e^{\frac{-\left|x-x_0\right|^2}{\hbar}}dx \right)\notag\\
&\leq& 4p_2| \nabla \Psi|_\infty^2e^{-4p_2\frac{\varrho^2}{\hbar^2}t }e^{\frac{-(\varrho-\epsilon)^2}{\hbar}}\int_\Omega \left|v(x,0)\right|^2dx\text{.}
\label{bubu}
\end{eqnarray}
Integrating (\ref{bubu}) over $(t;T)$ gives
\begin{eqnarray}
&~&\int_{\Omega \cap B(x_0,\varrho)} \left|\Psi(x) v(x,T)\right|^2 e^{\frac{-\left|x-x_0\right|^2}{\hbar}}dx\notag\\ &\leq& e^{4p_2\frac{\varrho^2}{\hbar^2}(T-t)}\int_{\Omega \cap B(x_0,\varrho)} \left|\Psi(x) v(x,t)\right|^2 e^{\frac{-\left|x-x_0\right|^2}{\hbar}}dx\notag\\
&~&+4p_2| \nabla \Psi|_\infty^2e^{4p_2\frac{\varrho^2}{\hbar^2}T}e^{\frac{-(\varrho-\epsilon)^2}{\hbar}}\int_\Omega \left|v(x,0)\right|^2dx \int_t^Te^{-4p_2\frac{\varrho^2}{\hbar^2}s}ds
\text{.}
\label{smile4}
\end{eqnarray}
It implies from the fact $\Psi(\cdot) =1$ in $B(x_0,\varrho-\epsilon)$ that
\begin{eqnarray}
\int_{\Omega \cap B(x_0,\varrho)} \left|\Psi(x) v(x,T)\right|^2 e^{\frac{-\left|x-x_0\right|^2}{\hbar}}dx 
&\geq& \int_{\Omega \cap B(x_0,\varrho - \epsilon)}\left|\Psi(x)v(x,T)\right|^2e^{\frac{-\left|x-x_0\right|^2}{\hbar}}dx\notag\\
&=& \int_{\Omega \cap B(x_0,\varrho - \epsilon)}\left|v(x,T)\right|^2e^{\frac{-\left|x-x_0\right|^2}{\hbar}}dx\notag\\
&\geq& 
\int_{\Omega \cap B(x_0,\varrho - 2\epsilon)}\left|v(x,T)\right|^2e^{\frac{-\left|x-x_0\right|^2}{\hbar}}dx\notag\\
&\geq& 
e^{\frac{-(\varrho-2\epsilon)^2}{\hbar}}\int_{\Omega \cap B(x_0,\varrho - 2\epsilon)}\left|v(x,T)\right|^2dx
\text{.}
\label{smile5}
\end{eqnarray}
Combining (\ref{smile4}), (\ref{smile5}) and the following estimate
\begin{equation}
\int_t^Te^{-4p_2\frac{\varrho^2}{\hbar^2}s}ds=\frac{1}{4p_2\frac{\varrho^2}{\hbar^2}} \left(e^{-4p_2\frac{\varrho^2}{\hbar^2}T}-e^{-4p_2\frac{\varrho^2}{\hbar^2}t}\right)\leq \frac{1}{4p_2\frac{\varrho^2}{\hbar^2}} e^{-4p_2\frac{\varrho^2}{\hbar^2}t}\leq \frac{1}{4p_2\varrho^2} e^{-4p_2\frac{\varrho^2}{\hbar^2}t}\end{equation} for $\hbar<1$, we obtain
\begin{eqnarray}
&~&\int_{\Omega \cap B(x_0,\varrho - 2\epsilon)}\left| v(x,T)\right|^2dx\notag\\
& \leq &e^{4p_2\frac{\varrho^2}{\hbar^2}(T-t)}e^{\frac{(\varrho-2\epsilon)^2}{\hbar}}\int_{\Omega \cap B(x_0,\varrho)} \left| v(x,t)\right|^2dx \notag \notag\\
& + &\frac{| \nabla \Psi|_\infty^2}{\varrho^2}e^{4p_2\frac{\varrho^2}{\hbar^2}(T-t)}e^{\frac{(\varrho-2\epsilon)^2}{\hbar}-\frac{(\varrho-\epsilon)^2}{\hbar}}\int_\Omega \left|v(x,0)\right|^2dx 
\text{.}
\end{eqnarray}
Let $\frac T2<T-\eta \hbar \leq t\leq T$, it yields
\begin{eqnarray}
\int_{\Omega \cap B(x_0,\varrho - 2\epsilon)}\left| v(x,T)\right|^2dx
& \leq &e^{4p_2\frac{\varrho^2\eta}{\hbar}}e^{\frac{(\varrho-2\epsilon)^2}{\hbar}}\int_{\Omega \cap B(x_0,\varrho)} \left| v(x,t)\right|^2dx \notag \notag\\
& + &\frac{| \nabla \Psi|_\infty^2}{\varrho^2}e^{4p_2\frac{\varrho^2\eta}{\hbar}}e^{\frac{(\varrho-2\epsilon)^2}{\hbar}-\frac{(\varrho-\epsilon)^2}{\hbar}}\int_\Omega \left|v(x,0)\right|^2dx 
\text{.}
\label{tam1}
\end{eqnarray}
We choose $\eta= \frac{\epsilon(2\varrho-3\epsilon)}{8p_2\varrho^2}$, that is $4p_2\varrho^2\eta=\frac 12\epsilon(-3\epsilon+2\varrho)$ in order to get
\[ 4p_2\varrho^2\eta + (\varrho-2\epsilon)^2 -(\varrho-\epsilon)^2 <0\]
Then, (\ref{tam1}) becomes
\begin{eqnarray}
\int_{\Omega \cap B(x_0,\varrho - 2\epsilon)}\left| v(x,T)\right|^2dx
& \leq &e^{\frac{(\varrho-\epsilon)^2+(\varrho-2\epsilon)^2}{2\hbar}}\int_{\Omega \cap B(x_0,\varrho)} \left| v(x,t)\right|^2dx \notag \notag\\
& + &\frac{| \nabla \Psi|_\infty^2}{\varrho^2}e^{\frac{(\varrho-2\epsilon)^2-(\varrho-\epsilon)^2}{2\hbar}}\int_\Omega \left|v(x,0)\right|^2dx 
\text{.}
\label{tam}
\end{eqnarray}
\textbf{Step 3:} Choose $\hbar$ for minimization problem.\\
Now, for the purpose of minimizing the right-hand side of inequality (\ref{tam}), we choose $\hbar$ such that
 \begin{equation*}
   \frac{| \nabla \Psi|_\infty^2}{\varrho^2}e^{\frac{(\varrho-2\epsilon)^2-(\varrho-\epsilon)^2}{2\hbar}}\int_\Omega \left|v(x,0)\right|^2dx \leq \frac{1}{2}\int_{\Omega \cap B(x_0,\varrho - 2\epsilon)}\left| v(x,T)\right|^2dx
 \end{equation*}
or 
 \begin{equation*}
 e^{\frac{(\varrho-2\epsilon)^2}{2\hbar}}  \leq e^{\frac{(\varrho-\epsilon)^2}{2\hbar}}\frac{\varrho^2}{2| \nabla \Psi|_\infty^2}\frac{\int_{\Omega \cap B(x_0,\varrho - 2\epsilon)}\left| v(x,T)\right|^2dx}{\int_\Omega \left|v(x,0)\right|^2dx}\text{.}
 \end{equation*}
 The choice of $\hbar$ also satisfies the condition that $\hbar <1$ and $\eta \hbar < \min\{1, \frac T2\}$. Such $\hbar$ exists by choosing
 \[ \hbar = \frac{\frac{\epsilon(2\varrho - 3\epsilon)}{2}}{\ln\left(\frac{2| \nabla \Psi|_\infty^2}{\varrho^2}\frac{\int_\Omega \left|v(x,0)\right|^2dx}{\int_{\Omega \cap B(x_0,\varrho - 2\epsilon)}\left| v(x,T)\right|^2dx}\right)+\frac{\epsilon(2\varrho - 3\epsilon)}{2}\left[1+\eta\left(1+\frac 2T\right)\right]}\text{.}\]
With this choice, it implies from (\ref{tam}) that
 \begin{equation*}
  \int_{\Omega \cap B(x_0,\varrho - 2\epsilon)}\left| v(x,T)\right|^2dx \leq \frac{\varrho^2}{| \nabla \Psi|_\infty^2}e^{\frac{(\varrho-2\epsilon)^2}{\hbar}}\frac{\int_{\Omega \cap B(x_0,\varrho - 2\epsilon)}\left| v(x,T)\right|^2dx}{\int_\Omega \left|v(x,0)\right|^2dx}\int_{\Omega \cap B(x_0,\varrho)} \left| v(x,t)\right|^2dx
    \text{.}\end{equation*}
This is equivalent to
 \begin{equation*}
\frac {\int_\Omega \left|v(x,0)\right|^2dx}{\int_{\Omega \cap B(x_0,\varrho)} \left| v(x,t)\right|^2dx
} \leq \frac{\varrho^2}{| \nabla \Psi|_\infty^2}e^{\frac{(\varrho-2\epsilon)^2}{\hbar}}\\
  \text{.}\end{equation*}
 This completes the proof with $E_1=\max\left\{\frac{\varrho^2}{| \nabla \Psi|_\infty^2};\frac{2(\varrho-2\epsilon)^2}{\epsilon(2\varrho - 3\epsilon)};1\right\}>1$ and
 \begin{eqnarray*}
 \frac 1\theta &=& \ln\left(\frac{2| \nabla \Psi|_\infty^2}{\varrho^2}\frac{2\int_\Omega \left|v(x,0)\right|^2dx}{\int_{\Omega \cap B(x_0,\varrho - 2\epsilon)}\left| v(x,T)\right|^2dx}\right)+\frac{\epsilon(2\varrho - 3\epsilon)}{2}\left[1+\eta\left(1+\frac 2T\right)\right]\\
 &=& \ln\left(\frac{2| \nabla \Psi|_\infty^2}{\varrho^2}\frac{2\int_\Omega \left|v(x,0)\right|^2dx}{\int_{\Omega \cap B(x_0,\varrho - 2\epsilon)}\left| v(x,T)\right|^2dx}e^{\frac{\epsilon(2\varrho - 3\epsilon)}{2}\left[1+\eta\left(1+\frac 2T\right)\right]}\right)
 \end{eqnarray*}
 Thus, there exists a constant $E_3>0$ such that
 \begin{equation*}  \frac 1\theta=\ln \left(E_3e^{\frac {E_3}T}\frac{\int_\Omega \left|v(x,0)\right|^2dx}{\int_{\Omega \cap B(x_0,\varrho - 2\epsilon)}\left| v(x,T)\right|^2dx}\right)
 \end{equation*}
 On the other hand, $\eta \hbar = \theta \frac{\epsilon^2(2\varrho-3\epsilon)^2}{16p_2\varrho^2}$. Hence, $E_2=\frac{\epsilon^2(2\varrho-3\epsilon)^2}{16p_2\varrho^2}$.
\end{proof}
 \subsection{Proof of Theorem \ref{convex}}
\begin{proof}
Let $x_0\in \Omega$, $\rho >0$ and $\ell >1$ such that $\ell \rho \leq \min\{\frac 12, \frac T4\}$, Corollary \ref{col1} gives
 \begin{equation}
\left|y'(t)+2N(t)y(t)\right| \leq \left(\frac{C_0}{T-t+\rho}+C_1 \right)y(t) 
\end{equation}
and \begin{equation}
N'(t)\leq\left(\frac{1+C_0}{T-t+\rho}+C_1\right)N(t)\text{.}
\end{equation}
Here
\begin{equation*}
y(t)=\frac{1}{(T-t+\rho)^{\frac{n}{2}}}\int_{\Omega} |v(x,t)|^2e^{\frac{-|x-x_0|^2}{4p(T)(T-t+\rho)}}dx \text{,}
\end{equation*}
\begin{equation*}
N(t)=\frac{p(t)}{(T-t+\rho)^{\frac{n}{2}}}\frac{\int_{\Omega} |\nabla v(x,t)|^2e^{\frac{-|x-x_0|^2}{4p(T)(T-t+\rho)}}dx}{y(t)}\text{,}
\end{equation*}
\begin{equation*}C_0=\frac{|p'|_\infty R^2 }{2p_1^2}~~\text{and}~~ C_1=\left(2+n\right)\frac{|p'|_\infty}{p_1}\text{.}
\end{equation*}
Now, thanks to Corollary \ref{col2}, one has
\begin{eqnarray}
y(T-\ell \rho)^{1+M_\ell } \leq e^{G_\ell } (1+2\ell )^{C_0(1+M_\ell )} y(T)y(T-2\ell \rho)^{M_\ell }
\label{bla}
\end{eqnarray}
where 
\begin{equation*}G_\ell  = C_1(1+M_\ell )\end{equation*}
and \begin{equation*}M_\ell \leq S_\ell = e^{C_1} \left\{
\begin{array}{ll}
\frac{\ln(1+\ell )}{\ln\frac{3}{2}}&\text{ if}~~ C_0=0\text{,}\\
\frac{(1+\ell )^{C_0}}{1-\left(\frac{2}{3}\right)^{C_0}} &\text{if}~~ C_0>0\text{.}
\end{array} \right. \end{equation*}
Then (\ref{bla}) is equivalent to
\begin{eqnarray}
 &~&\left(\int_{\Omega} |v(x,T-\ell \rho)|^2dx\right)^{1+M_\ell}\notag\\
&\leq& \frac{\left((1+\ell)\rho\right)^{(1+M_\ell)\frac n2}}{\rho^{\frac n2}\left((1+2\ell)\rho\right)^{M_\ell\frac n2}}e^{\frac{R^2(1+M_\ell)}{4p(T)(1+\ell)\rho}}e^{C_1(1+M_\ell)}(1+2\ell)^{C_0(1+M_\ell)}\notag\\
&~&\times \left(\int_{\Omega} |v(x,t)|^2 e^{\frac{-|x-x_0|^2}{4\rho p(T)}}dx\right)\left(\int_{\Omega} |v(x,T-2\ell \rho)|^2e^{\frac{-|x-x_0|^2}{4 (1+\ell)\rho p(T)}}dx\right)^{M_\ell}\notag\\
&\leq& 2^{C_0(1+M_\ell)} (1+\ell)^{\frac{n}{2}+C_0(1+M_\ell)}e^{(1+M_\ell)C_1}e^{\frac{R^2(1+M_\ell)}{4p(T)(1+\ell)\rho}}\notag\\
&~&\times\left(\int_{\Omega} |v(x,t)|^2 e^{\frac{-|x-x_0|^2}{4\rho p(T)}}dx\right)\left(\int_{\Omega} |v(x,0)|^2dx\right)^{M_\ell} \text{.}
\label{111}
\end{eqnarray}
On the other hand, for  $0<r<R$ such that $B(x_0,r)\subset \omega$, one has
\begin{eqnarray}
\int_{\Omega} |v(x,T)|^2 e^{\frac{-|x-x_0|^2}{4\rho p(T)}}dx
&\leq& \int_{ B(x_0,r)} |v(x,T)|^2dx + \int_{\Omega \cap\{x;|x-x_0|\geq r\}} |v(x,T)|^2 e^{\frac{-|x-x_0|^2}{4\rho p(T)}}dx\notag\\
&\leq& \int_{B(x_0,r)} |v(x,T)|^2dx +e^{\frac{-r^2}{4\rho p(T)}} \int_{\Omega} |v(x,T)|^2dx \notag\\
&\leq& \int_\omega |v(x,T)|^2dx +e^{\frac{-r^2}{4\rho p(T)}}\int_{\Omega} |v(x,0)|^2dx 
\text{.}
\label{112}
\end{eqnarray}
Moreover, we also have
\begin{equation}
\int_{\Omega} |v(x,T-\ell \rho)|^2dx\geq\int_{\Omega} |v(x,T)|^2 e^{\frac{-|x-x_0|^2}{4\rho p(T)}}dx
\label{113}
\end{equation}
Combining (\ref{111}), (\ref{112}) and (\ref{113}), it yields
\begin{eqnarray}
&~&\left(\int_\Omega |v(x,T)|^2dx\right)^{1+M_\ell}\notag\\
&\leq& 2^{C_0(1+M_\ell)} (1+\ell)^{\frac{n}{2}+C_0(1+M_\ell)}e^{(1+M_\ell)C_1}e^{\frac{R^2(1+M_\ell)}{4(1+\ell)\rho p(T)}}\notag\\
&\times&\left[\int_{\omega} |v(x,T)|^2dx\left(\int_\Omega |v(x,0)|^2dx\right)^{M_\ell}+e^{-\frac{r^2}{4\rho p(T)}}\left(\int_\Omega |v(x,0)|^2dx\right)^{1+M_\ell}\right]\text{.}\notag\\
\label{min}
\end{eqnarray}
With the notice that $M_\ell \leq S_\ell$, we can write (\ref{min}) as below, thanks to the energy estimate $\int_\Omega |v(x,0)|^2dx \geq \int_\Omega |v(x,T)|^2dx$
\begin{eqnarray}
&~&\left(\int_\Omega |v(x,T)|^2dx\right)^{1+S_\ell}\notag\\
&\leq& 2^{C_0(1+S_\ell)} (1+\ell)^{\frac{n}{2}+C_0(1+S_\ell)}e^{(1+S_\ell)C_1}e^{\frac{R^2(1+S_\ell)}{4(1+\ell)\rho p(T)}}\notag\\
&\times&\left[\int_{\omega} |v(x,T)|^2dx\left(\int_\Omega |v(x,0)|^2dx\right)^{S_\ell}+e^{-\frac{r^2}{4\rho p(T)}}\left(\int_\Omega |v(x,0)|^2dx\right)^{1+S_\ell}\right]\text{.}\notag\\
\label{min2}
\end{eqnarray}
Now, in order to minimize the right-hand side of inequality (\ref{min2}), we will choose $\ell>1$ as
\begin{equation*}
\ell =\left\{ \begin{array}{ll}
\left(\frac{2^{2+\xi} R^2e^{C_1}}{\xi \ln \frac{3}{2}r^2}\right)^{\frac{1}{1-\xi}}-1 ~~\forall \xi \in (0,1)~~&\text{if}~~~C_0 = 0\text{,}\\
\left(\frac{4R^2e^{C_1}}{r^2\left(1-\left(\frac{2}{3}\right)^{C_0}\right)}\right)^{\frac{1}{1-C_0}}-1 ~~&\text{if}~~~C_0 > 0 \text{.}
\end{array}
\right.
\end{equation*}
The assumption (i) follows $C_0 <1$. Hence, such choice of $\ell$ provides us
\begin{equation*}
\frac{R^2(1+S_\ell)}{4(1+\ell)\rho p(T)}\leq \frac{r^2}{8\rho p(T)}\text{.}
\end{equation*}
Thus, (\ref{min2}) becomes:
 \begin{eqnarray}
&~&\left(\int_\Omega |v(x,T)|^2dx\right)^{1+S_\ell}\notag\\
&\leq& 2^{C_0(1+S_\ell)} (1+\ell)^{\frac{n}{2}+C_0(1+S_\ell)}e^{(1+S_\ell)C_1}e^{\frac{r^2}{8p(T)h}}\int_{\omega} |v(x,T)|^2dx\left(\int_\Omega |v(x,0)|^2dx\right)^{S_\ell}\notag\\
&+&2^{C_0(1+S_\ell)} (1+\ell)^{\frac{n}{2}+C_0(1+S_\ell)}e^{(1+S_\ell)C_1}e^{\frac{-r^2}{8p(T)h}}\left(\int_\Omega |v(x,0)|^2dx\right)^{1+S_\ell}
\text{.}
\label{step71}
\end{eqnarray}
This estimate is true for any $\rho >0$ satisfying $\rho <\frac{1}{\ell} \min\{\frac{1}{2}; \frac{T}{4}\}$. For $\rho \geq\frac{1}{\ell} \min\{\frac{1}{2}; \frac{T}{4}\}$ which implies $\frac{r^2}{8\rho p(T)}\leq \frac{r^2\ell} {4p(T)}\left(1+\frac 2T\right)$, we can get the following estimate be true for any $\rho >0$.
\begin{eqnarray}
&~&\left(\int_\Omega |v(x,T)|^2dx\right)^{1+S_\ell}\notag\\
&\leq& 2^{C_0(1+S_\ell)} (1+\ell)^{\frac{n}{2}+C_0(1+S_\ell)}e^{(1+S_\ell)C_1}e^{\frac{r^2}{8\rho p(T)}}\left(\int_{\omega} |v(x,T)|^2dx\right)\left(\int_\Omega |v(x,0)|^2dx\right)^{S_\ell} \notag\\
&+&2^{C_0(1+S_\ell)} (1+\ell)^{\frac{n}{2}+C_0(1+S_\ell)}e^{(1+S_\ell)C_1}e^{\frac{r^2\ell}{4p(T)}\left(1+\frac 2T\right)}e^{\frac{-r^2}{8\rho p(T)}}\left(\int_\Omega |v(x,0)|^2dx\right)^{1+S_\ell} \text{.}
\end{eqnarray}
Now, we choose $\rho$ such that
\begin{eqnarray*}
&~&2^{C_0(1+S_\ell)} (1+\ell)^{\frac{n}{2}+C_0(1+S_\ell)}e^{(1+S_\ell)C_1}e^{\frac{r^2\ell}{4p(T)}\left(1+\frac 2T\right)}e^{\frac{-r^2}{8\rho p(T)}}\left(\int_\Omega |v(x,0)|^2dx\right)^{1+S_\ell}\\
&=&\frac{1}{2}\left(\int_\Omega |v(x,T)|^2dx\right)^{1+S_\ell}
\end{eqnarray*}
that is
\begin{equation*}
e^{\frac{r^2}{8\rho p(T)}}=22^{C_0(1+S_\ell)} (1+\ell)^{\frac{n}{2}+C_0(1+S_\ell)}e^{(1+S_\ell)C_1}e^{\frac{r^2\ell}{4p(T)}\left(1+\frac 2T\right)}\left(\frac{\int_\Omega |v(x,0)|^2dx}{\int_\Omega |v(x,T)|^2dx}\right)^{1+S_\ell}\text{.}
\end{equation*}
Therefore, we get\\
\begin{eqnarray*}
\left(\int_\Omega |v(x,T)|^2dx\right)^{2(1+S_\ell)}
&\leq& 4.4^{C_0(1+S_\ell)} (1+\ell)^{n+2C_0(1+S_\ell)}e^{2(1+S_\ell)C_1}e^{\frac{r^2\ell}{4p(T)}\left(1+\frac 2T\right)}\\
&\times&\left(\int_{\omega} |v(x,T)|^2dx\right)\left(\int_\Omega |v(x,0)|^2dx\right)^{1+2S_\ell} \text{.}
\end{eqnarray*}
Thus, we can state
\begin{eqnarray}
&~&\int_\Omega |v(x,T)|^2dx\notag\\
&\leq&  \left(4^{1+C_0(1+S_\ell)} (1+\ell)^{n+2C_0(1+S_\ell)}e^{2(1+S_\ell)C_1}e^{\frac{r^2\ell}{4p(T)}\left(1+\frac 2T\right)}\int_{\omega} |v(x,T)|^2dx\right)^{\frac{1}{2(1+S_\ell)}}\notag\\
&~& \times \left(\int_\Omega |v(x,0)|^2dx\right)^{\frac{1+2S_\ell}{2(1+S_\ell)}}\text{.}
\end{eqnarray}
This completes the proof.
 \end{proof}
\subsection{Proof of Theorem \ref{nonconvex}}
Let us move to the proof of Theorem \ref{nonconvex} with the structure as: Thanks to Preliminary results in the previous subsection, we can get a H\"older type estimate in Step 1. Furthermore, thanks to the technical Lemma \ref{lemma3}, we get an estimate for the term containing $F(t)$, which is presented in Step 2. Step 3 will make appear a small ball, which is related to the presence of $\omega$ later, by using a splitting technique. Next, dealing with a minimization problem, Step 4 provides us a localized observation estimate. Due to propagation of smallness by constructing the sequence of balls, $\omega$ will appear in Step 5. Lastly, in Step 6, by using an adequate covering of $\Omega$ with a finite number of balls, we will get the desired result.
\begin{proof}
\textbf{Step 1:} Get H\"{o}lder type inequality.\\
Let $x_0 \in \Omega$ and $R$ be small enough such that $\Omega \cap B(x_0;R)$ is star-shaped with respect to $x_0$. Such choice of $x_0$ and $R$ will be mentioned in Step 6. Let $0<\epsilon<\frac{R}{4}$.\\ Define $\psi \in C_0^2(B(x_0;R))$ satisfying $\psi = 1$ in $B(x_0;R-\epsilon)$ and $0<\psi(t)<1 ~~ \forall t\in B(x_0;R)$.
Then $\psi v \in H^1((0,T);H_0^1(\Omega \cap B(x_0;R))$. For any $\rho >0$ and $\ell >1$ such that $\ell \rho \leq \min\{\frac 12, \frac T4\}$, Corollary \ref{col1} gives
 \begin{equation}
\left|y'(t)+2N(t)y(t)\right| \leq \left(\frac{C_0}{T-t+\rho}+C_1 \right)y(t) +2\int_\vartheta |w(x,t)z(x,t)|e^{\phi(x,t)}dx
\label{hehe}
\end{equation}
and \begin{equation}
N'(t)\leq\left(\frac{1+C_0}{T-t+\rho}+C_1\right)N(t)+\frac{1}{2}\frac{ \int_{\vartheta} |w(x,t)|^2e^{\phi(x,t)}dx}{y(t)}\text{.}
\end{equation}
Here
\[\vartheta = \Omega \cap B(x_0;R)\text{,}\]
\[z=\psi v\text{,}\]
\[w:=\partial_t z -\Delta z\text{,}\]
\begin{equation*}
y(t)=\frac{1}{(T-t+\rho)^{\frac{n}{2}}}\int_{\vartheta} |z(x,t)|^2e^{\frac{-|x-x_0|^2}{4p(T)(T-t+\rho)}}dx \text{,}
\end{equation*}
\begin{equation*}
N(t)=\frac{p(t)}{(T-t+\rho)^{\frac{n}{2}}}\frac{\int_{\vartheta} |\nabla z(x,t)|^2e^{\frac{-|x-x_0|^2}{4p(T)(T-t+\rho)}}dx}{y(t)}\text{,}
\end{equation*}
\begin{equation*}C_0=\frac{|p'|_\infty R^2 }{2p_1^2}~~\text{and}~~ C_1=\left(2+n\right)\frac{|p'|_\infty}{p_1}\text{.}
\end{equation*}
Thanks to the Cauchy-Schwarz inequality and the fact that $2ab \leq a^2+b^2~~\forall a,b$, one gets
\begin{eqnarray}
&~&2\int_\vartheta |w(x,t)z(x,t)|e^{\phi(x,t)}dx\notag\\
&\leq& 2\left(\int_\vartheta |w(x,t)|^2e^{\phi(x,t)}dx\right)^{\frac 12}\left(\int_\vartheta |z(x,t)|^2e^{\phi(x,t)}dx\right)^{\frac 12}\notag\\
&\leq& \int_\vartheta |w(x,t)|^2e^{\phi(x,t)}dx + \int_\vartheta |z(x,t)|^2e^{\phi(x,t)}dx\text{,}
\end{eqnarray}
Thus, we can write (\ref{hehe}) as below
\begin{equation}
\left|y'(t)+2N(t)y(t)\right| \leq \left(\frac{C_0}{T-t+\rho}+C_1+1+F(t) \right)y(t)
\end{equation}
with \[F(t):=\frac{ \int_{\vartheta} |w(x,t)|^2e^{\phi(x,t)}dx}{y(t)}\text{.}\]
Now, thanks to Corollary \ref{col2}, one has
\begin{eqnarray}
y(T-\ell \rho)^{1+M_\ell } \leq e^{G_\ell } (1+2\ell )^{C_0(1+M_\ell )} y(T)y(T-2\ell \rho)^{M_\ell }
\label{step11}
\end{eqnarray}
where 
\begin{equation*}G_\ell  = (1+M_\ell )\left(C_1 +1 + 2\int_{T-2\ell \rho}^{T}F(s)ds\right)\end{equation*}
and
\begin{equation*}M_\ell \leq S_\ell = e^{C_1+1}\left\{
\begin{array}{ll}
\frac{\ln(1+\ell )}{\ln\frac{3}{2}}&\text{ if}~~ C_0=0 \text{,}\\
\frac{(1+\ell )^{C_0}}{1-\left(\frac{2}{3}\right)^{C_0}} &\text{if}~~ C_0>0\text{.}
\end{array} \right. \end{equation*}

 \textbf{Step 2:} Estimate $\int_{T-2\ell \rho}^{T} F(s)ds$. \\
Remind that
\begin{eqnarray}
 F(s)&=& \frac{\int_{\vartheta} |w(x,s)|^2e^{\phi(x,s)}dx}{\int_{\vartheta} |z(x,s)|^2e^{\phi(x,s)}dx}
 \label{F}
  \end{eqnarray}
 with 
 \begin{eqnarray}
 w(x,s)&=& \partial_t z(x,s) -p(s)\Delta z(x,s)\notag\\
 &=& -p(s)v(x,s)\Delta \psi(x) -2p(s) \nabla v(x,s) \nabla \psi(x)\text{.}
 \end{eqnarray}
 Note that $\nabla \psi = \Delta \psi =0$ in $\Omega \cap B(x_0,R-\epsilon)$, so
 \begin{eqnarray}
&~&\int_{\vartheta} |w(x,s)|^2e^{\phi(x,s)}dx\notag\\
&=& p(s)^2\int_{\vartheta} \left(v(x,s)\Delta \psi(x) +2\nabla v(x,s)) \nabla \psi(x)\right)^2e^{\phi(x,s)}dx\notag\\
&=& p(s)^2\int_{\Omega \cap \{x; |x-x_0|\geq R-\epsilon\}}\left(v(x,s)\Delta \psi(x) +2\nabla v(x,s)) \nabla \psi(x)\right)^2e^{\phi(x,s)}dx\text{.}
 \end{eqnarray}
It implies from the fact that $e^{\phi(x,s)}=\frac{1}{(T-s+\rho)^{\frac{n}{2}}}e^{\frac{-|x-x_0|^2}{4p(T)(T-s+\rho)}}$ and $(a+b)^2\leq 2(a^2+b^2)~~\forall a,b$ that
 \begin{eqnarray}
&~&\int_{\vartheta} |w(x,s)|^2e^{\phi(x,s)}dx\notag\\
&\leq&\frac{2p(s)^2}{(T-s+\rho)^{\frac{n}{2}}}e^{\frac{-(R-\epsilon)^2}{4p(T)(T-s+\rho)}}\left(|\Delta \psi|_\infty^2\int_\Omega |v(x,s)|^2dx+4|\nabla \psi|_\infty^2\int_\Omega |\nabla v(x,s)|^2dx\right)\text{.}\notag\\
 \end{eqnarray}
 Moreover, thanks to the following energy estimate 
  \begin{equation*}
 \int_\Omega |\nabla v(x,s)|^2dx \leq \frac{1}{2p_1s}\int_\Omega |v(x,0)|^2dx~~\forall s>0\text{,}
 \end{equation*}
we obtain: 
\begin{equation}
\int_{\vartheta} |w(x,s)|^2e^{\phi(x,s)}dx\leq \frac{C_2}{(T-s+\rho)^{\frac{n}{2}}}e^{\frac{-(R-\epsilon)^2}{4p(T)(T-s+\rho)}}\left(1+\frac{1}{2p_1s}\right)\int_\Omega |v(x,0)|^2dx
\label{f1}
\end{equation}
where $C_2 = 2p_2^2 \max \{|\Delta \psi|_\infty^2, 4|\nabla \psi|_\infty^2\}$. On the other hand, due to $\psi =1$ on $B(x_0,R-\epsilon)$, one gets
\begin{eqnarray}
\int_\vartheta |z(x,s)|^2e^{\phi(x,s)}dx &=& \int_{\Omega\cap B(x_0,R-\epsilon)} |v(x,s)|^2e^{\phi(x,s)}dx\notag\\
&\geq& \frac{1}{(T-s+\rho)^{\frac{n}{2}}}e^{\frac{-(R-2\epsilon)^2}{4p(T)(T-s+\rho)}}\int_{\Omega\cap B(x_0,R-2\epsilon)} |v(x,s)|^2dx\text{.}
\label{f2}
\end{eqnarray}
Combining (\ref{F}), (\ref{f1}) and (\ref{f2}) gives
\begin{equation}
F(s) \leq C_2e^{\frac{-\epsilon(2R-3\epsilon)}{12\ell \rho p(T)}}\left(1+\frac{1}{2p_1s}\right)\frac{\int_\Omega |v(x,0)|^2dx}{\int_{\Omega\cap B(x_0,R-2\epsilon)} |v(x,s)|^2dx}
\label{f3}
\end{equation}
with $s\in[T-2\ell \rho,T]$.
Now, apply Lemma \ref{lemma3} with $\varrho = R-2\epsilon$, under condition $2\ell \rho \leq E_2\theta$ for some $E_2 >0$ depending on $R$ and $\epsilon$,
  there exists a constant $E_1=E_1(R,\epsilon)>1$ such that
 \begin{equation}
  \frac{\int_\Omega \left|v(x,0)\right|^2dx}{\int_{\Omega \cap B(x_0,R-2\epsilon)} \left|v(x,s)\right|^2dx} \leq E_1e^{\frac{E_1}{\theta}} ~~\forall T-\ell \rho\leq s \leq T\text{.}
\label{f4}
  \end{equation}
Thus, from (\ref{f3}) and (\ref{f4}), one has
\begin{eqnarray*}
\int_{T-2\ell \rho}^{T} F(s)ds \leq E_1e^{\frac{E_1}{\theta} }C_2e^{\frac{-\epsilon(2R-3\epsilon)}{12\ell \rho p_1}}\int_{T-2\ell \rho}^{T} \left(1+\frac{1}{2p_1s}\right)ds\text{.}
\end{eqnarray*}
In order to get $\frac{E_1}{\theta}-\frac{\epsilon(2R-3\epsilon)}{12\ell \rho p_1}<0$, we take $2\ell \rho \leq cE_2\theta$ with $c=\min\left\{\frac{\epsilon(2R-3\epsilon)}{6p_1E_1E_2};1\right\}$. On the other hand, due to $\int_{T-2\ell \rho}^{T} \left(1+\frac{1}{2p_1s}\right)ds = 2\ell \rho + \frac{1}{2p_1}\ln \frac{T}{T-2\ell \rho} \leq 1+\frac{1}{2p_1}\ln 2$,  there exists a constant $C_3>0$ which does not depend on $\ell,\rho$ and $T$ such that 
\begin{equation*}
\int_{T-2\ell \rho}^{T} F(s)ds \leq C_3 \text{.}
\end{equation*}

\textbf{Step 3:} Make appear a small ball.\\
Remind that
\begin{equation*}
y(t)=\int_\vartheta |z(x,t)|^2 \frac{1}{(T-t+\rho)^{\frac{n}{2}}}e^{\frac{-|x-x_0|^2}{4p(T)(T-t+\rho)}}\text{.}
\end{equation*}
From (\ref{step11}) in Step 1, we have
\begin{eqnarray}
 &~&\left(\int_{\vartheta} |z(x,T-\ell \rho)|^2dx\right)^{1+M_\ell}\notag\\
&\leq& (1+2\ell)^{\frac{n}{2}+C_0(1+M_\ell)}e^{(1+M_\ell)(C_1+1+2C_3)}e^{\frac{R^2(1+M_\ell)}{4(1+\ell)\rho p(T)}}\notag\\
&~&\left(\int_{\vartheta} |z(x,t)|^2 e^{\frac{-|x-x_0|^2}{4\rho p(T)}}dx\right)\left(\int_{\vartheta} |z(x,T-2\ell \rho)|^2e^{\frac{-|x-x_0|^2}{4(1+2\ell)\rho p(T)}}dx\right)^{M_\ell} \text{.}
\end{eqnarray}
 From the fact that $|z| \leq |v|$ in $\Omega \cap B(x_0,R)$ and $\int_\Omega |v(x,T-2\ell \rho)|^2dx\leq \int_\Omega |v(x,0)|^2dx$, we can write
\begin{eqnarray}
 &~&\left(\int_{\vartheta} |z(x,T-\ell \rho)|^2dx\right)^{1+M_\ell}\notag\\
&\leq& (1+2\ell)^{\frac{n}{2}+C_0(1+M_\ell)}e^{(1+M_\ell)(C_1+1+2C_3)}e^{\frac{R^2(1+M_\ell)}{4(1+\ell)\rho p(T)}}\left(\int_{\vartheta} |v(x,t)|^2 e^{\frac{-|x-x_0|^2}{4\rho p(T)}}dx\right)\left(\int_{\Omega}|v(x,0)|^2dx\right)^{M_\ell}\text{.} \notag\\
\label{44}
\end{eqnarray}
On the other hand, thanks to $\psi =1$ on $B(x_0, R-\epsilon)$, we obtain
\begin{eqnarray}
\int_\vartheta |z(x,T-\ell \rho)|^2dx \geq\int_{\Omega \cap B(x_0,R-2\epsilon)}|v(x,T-\ell \rho)|^2dx\text{.}
\label{45}
\end{eqnarray}
From inequality (\ref{f4}), one gets
\begin{equation}
\int_{\Omega \cap B(x_0,R-2\epsilon)}|v(x,T-\ell \rho)|^2dx \geq \frac{1}{E_1e^{\frac{E_1}{\theta}}}\int_\Omega |v(x,0)|^2dx \geq \frac{1}{E_1e^{\frac{E_1}{\theta}}}\int_\Omega |v(x,T)|^2dx\text{.}
\label{46}
\end{equation}
Combining (\ref{44}), (\ref{45}) and (\ref{46}), it yields
\begin{eqnarray}
\left(\int_\Omega |v(x,T)|^2dx\right)^{1+M_\ell}&\leq& (E_1e^{\frac{E_1}{\theta}})^{1+M_\ell}(1+2\ell)^{\frac{n}{2}+C_0(1+M_\ell)}e^{(1+M_\ell)(C_1+1+2C_3)}e^{\frac{R^2(1+M_\ell)}{4(1+\ell)\rho p(T)}}\notag\\
&~&\left(\int_{\vartheta} |v(x,t)|^2 e^{\frac{-|x-x_0|^2}{4\rho p(T)}}dx\right)\left(\int_{\Omega}|v(x,0)|^2dx\right)^{M_\ell}\text{.}
\end{eqnarray}
Furthermore, for  $0<r\leq \frac{R}{2}$ such that $B(x_0,r)\subset \Omega$, one has
\begin{eqnarray}
\int_{\vartheta} |v(x,T)|^2 e^{\frac{-|x-x_0|^2}{4\rho p(T)}}dx
&\leq& \int_{ B(x_0,r)} |v(x,T)|^2dx + \int_{\Omega \cap\{x;|x-x_0|\geq r\}} |v(x,T)|^2 e^{\frac{-|x-x_0|^2}{4\rho p(T)}}dx\notag\\
&\leq& \int_{B(x_0,r)} |v(x,T)|^2dx +e^{\frac{-r^2}{4\rho p(T)}}\int_{\Omega} |v(x,0)|^2dx 
\text{.}
\end{eqnarray}
Thus,  we obtain that\\
\begin{eqnarray}
&~&\left(\int_\Omega |v(x,T)|^2dx\right)^{1+M_\ell}\notag\\
&\leq& (E_1e^{\frac{E_1}{\theta}})^{1+M_\ell}(1+2\ell)^{\frac{n}{2}+C_0(1+M_\ell)}e^{(1+M_\ell)(C_1+1+2C_3)}\notag\\
&~&\left[e^{\frac{R^2(1+M_\ell)}{4(1+\ell)\rho p(T)}}\int_{\vartheta} |v(x,T)|^2dx\left(\int_\Omega |v(x,0)|^2dx\right)^{M_\ell}+e^{\frac{R^2(1+M_\ell)}{4(1+\ell)\rho p(T)}}e^{-\frac{r^2}{4\rho p(T)}}\left(\int_\Omega |v(x,0)|^2dx\right)^{1+M_\ell}\right]
\text{.}\notag\\
\label{step6}
\end{eqnarray}
Notice that $M_\ell \leq S_\ell$, therefore we can write (\ref{step6}) as below
\begin{eqnarray}
&~&\left(\int_\Omega |v(x,T)|^2dx\right)^{1+S_\ell}\notag\\
&\leq& (E_1e^{\frac{E_1}{\theta}})^{1+S_\ell}(1+2\ell)^{\frac{n}{2}+C_0(1+S_\ell)}e^{(1+S_\ell)(C_1+1+2C_3)}\notag\\
&~&\left[e^{\frac{R^2(1+S_\ell)}{4(1+\ell)\rho p(T)}}\int_{\vartheta} |v(x,T)|^2dx\left(\int_\Omega |v(x,0)|^2dx\right)^{S_\ell}+e^{\frac{R^2(1+S_\ell)}{4(1+\ell)\rho p(T)}}e^{-\frac{r^2}{4\rho p(T)}}\left(\int_\Omega |v(x,0)|^2dx\right)^{1+S_\ell}\right]
\text{.}\notag\\
\label{fun}
\end{eqnarray}

\textbf{Step 4:} Choose suitable $\ell,\rho$ and solve a minimization problem.\\
With $R$ is small enough such that $C_0:=\frac{R^2|p'|_\infty}{2p_1^2}<1$, we will minimize the right-hand side of the estimate (\ref{fun}) by choosing $\ell$ as below
\begin{equation*}
\ell =\left\{ \begin{array}{ll}
\left(\frac{2^{2+\xi} R^2e^{C_1+1}}{\xi \ln \frac{3}{2}r^2}\right)^{\frac{1}{1-\xi}}-1 ~~\forall \xi \in (0,1)~~&\text{if}~~~C_0 = 0; \\
\left(\frac{4R^2e^{C_1+1}}{r^2\left(1-\left(\frac{2}{3}\right)^{C_0}\right)}\right)^{\frac{1}{1-C_0}}-1 ~~&\text{if}~~~C_0> 0
\end{array}
\right.
\end{equation*}
in order to get
\begin{equation*}
\frac{R^2(1+S_\ell)}{4(1+\ell)\rho p(T)}\leq \frac{r^2}{8\rho p(T)}\text{.}
\end{equation*}
With this choice of $\ell$ and the fact that $E_1>1$, we can conclude from (\ref{fun}) that: there exists a constant $C_4>1$ not depending on $\rho$ such that
 \begin{eqnarray*}
&~&\left(\int_\Omega |v(x,T)|^2dx\right)^{1+S_\ell}\notag\\
&\leq&C_4e^{\frac{C_4}{\theta}}e^{\frac{r^2}{8p(T)\rho}}\int_{\vartheta} |v(x,T)|^2dx\left(\int_\Omega |v(x,0)|^2dx\right)^{S_\ell}+C_4e^{\frac{C_4}{\theta}}e^{\frac{-r^2}{8p(T)\rho}}\left(\int_\Omega |v(x,0)|^2dx\right)^{1+S_\ell}
\text{.}
\end{eqnarray*}
This estimate is true for any $\rho >0$ satisfying $\rho \leq\frac{1}{\ell}\min\{\frac{1}{2}; \frac{T}{4}; cE_2\theta\}$. For $\rho >\frac{1}{\ell}\min\{\frac{1}{2}; \frac{T}{4}; cE_2\theta\}$ which implies $\frac{r^2}{8\rho p(T)}< \frac{r^2\ell}{8p(T)}\left(2+\frac{4}{T}+\frac{1}{cE_2\theta}\right)$, we can get the following estimate be true for any $\rho >0$.
\begin{eqnarray*}
\left(\int_\Omega |v(x,T)|^2dx\right)^{1+S_\ell}
&\leq& C_4e^{\frac{C_4}{\theta}}e^{\frac{r^2}{8\rho p(T)}}\left(\int_{B(x_0,r)} |v(x,T)|^2dx\right)\left(\int_\Omega |v(x,0)|^2dx\right)^{S_\ell} \notag\\
&~&+ C_4e^{\frac{C_4}{\theta}}e^{\frac{r^2\ell}{8p(T)}\left(2+\frac{4}{T}+\frac{1}{cE_2\theta}\right)}e^{\frac{-r^2}{8\rho p(T)}}\left(\int_\Omega |v(x,0)|^2dx\right)^{1+S_\ell} \text{.}
\end{eqnarray*}
Now, we choose $\rho$ such as
\begin{equation*}C_4e^{\frac{C_4}{\theta}}e^{\frac{r^2\ell}{8p(T)}\left(2+\frac{4}{T}+\frac{1}{cE_2\theta}\right)}e^{\frac{-r^2}{8\rho p(T)}}\left(\int_\Omega |v(x,0)|^2dx\right)^{1+S_\ell}=\frac{1}{2}\left(\int_\Omega |v(x,T)|^2dx\right)^{1+S_\ell} \text{;}
\end{equation*}
that is
\begin{equation*}
e^{\frac{r^2}{8\rho p(T)}}=2C_4e^{\frac{C_4}{\theta}}e^{\frac{r^2\ell}{8p(T)}\left(2+\frac{4}{T}+\frac{1}{cE_2\theta}\right)}\left(\frac{\int_\Omega |v(x,0)|^2dx}{\int_\Omega |v(x,T)|^2dx}\right)^{1+S_\ell} \text{;}
\end{equation*}
in order to get\\
\begin{eqnarray*}
&~&\left(\int_\Omega |v(x,T)|^2dx\right)^{2(1+S_\ell)}\notag\\
&\leq& 4 C_4^2e^{\frac{2C_4}{\theta}}e^{\frac{r^2\ell}{8p(T)}\left(2+\frac{4}{T}+\frac{1}{cE_2\theta}\right)}\left(\int_{B(x_0,r)} |v(x,T)|^2dx\right)\left(\int_\Omega |v(x,0)|^2dx\right)^{1+2S_\ell} \text{.}
\end{eqnarray*}
Thus,  there exists a constant $C_5$ not depending on $T$ and $\theta$ such that
\begin{equation*}
\int_\Omega |v(x,T)|^2dx
\leq  C_5e^{\frac{C_5}{\theta}}e^{\frac{C_5}{T}}\left(\int_{B(x_0,r)} |v(x,T)|^2dx\right)^{\frac{1}{2(1+S_\ell)}}\left(\int_\Omega |v(x,0)|^2dx\right)^{\frac{1+2S_\ell}{2(1+S_\ell)}}\text{.}
\end{equation*}
On the other hand, Lemma \ref{lemma3} says that there exists $E_3 > 0$ satisfying the following estimate
\begin{equation*}
e^{\frac{1}{\theta}} = E_3e^{\frac{E_3}{T}}\frac{\int_\Omega |v(x,0)|^2dx}{\int_{\Omega\cap B(x_0,R-4\epsilon)}|v(x,T)|^2dx}\text{.}
\end{equation*}
Hence, the following estimate holds
\begin{eqnarray*}
\int_\Omega |v(x,T)|^2dx
&\leq&  C_5E_3^{C_5}e^{\frac{C_5(E_3+1)}{T}}\left(\frac{\int_\Omega |v(x,0)|^2dx}{\int_{\Omega\cap B(x_0,R-4\epsilon)}|v(x,T)|^2dx}\right)^{C_5}\notag\\
&~&\times \left(\int_{B(x_0,r)} |v(x,T)|^2dx\right)^{\frac{1}{2(1+S_\ell)}}\left(\int_\Omega |v(x,0)|^2dx\right)^{\frac{1+2S_\ell}{2(1+S_\ell)}}\text{.}
\end{eqnarray*}
Using the fact that $\int_{\Omega\cap B(x_0,R-4\epsilon)}|v(x,T)|^2dx\leq \int_{\Omega}|v(x,T)|^2dx$, we obtain
\begin{eqnarray*}
&~&\left(\int_{\Omega\cap B(x_0,R-4\epsilon)} |v(x,T)|^2dx\right)^{1+C_5}\notag\\
&\leq& C_5E_3^{C_5}e^{\frac{C_5(E_3+1)}{T}}\left(\int_{B(x_0,r)} |v(x,T)|^2dx\right)^{\frac{1}{2(1+S_\ell)}}\left(\int_\Omega |v(x,0)|^2dx\right)^{C_5+\frac{1+2S_\ell}{2(1+S_\ell)}}\text{.}
\end{eqnarray*}
Thus, there exist constants $\kappa>0$ and $\sigma \in (0,1)$ satisfying the following estimate
\begin{equation}
\int_ {\Omega\cap B(x_0,R-4\epsilon)}|v(x,T)|^2dx
\leq  \kappa e^{\frac{\kappa }{T}}\left(\int_{B(x_0,r)} |v(x,T)|^2dx\right)^{\sigma}\left(\int_{\Omega} |v(x,0)|^2dx\right)^{1-\sigma}\text{.}
\label{step7}
\end{equation}

\textbf{Step 5} Make appear $\omega$ by propagation of smallness.\\
Let $r>0$ be small enough and $x_j\in \Omega (j=1,2,...,m) (m\in \mathbb{N})$, we can construct a sequence of balls $\{B(x_j,r)\}_{j\in \overline{1,m}}$ such that the following inclusions hold
\begin{enumerate}
\item $B(x_m,r) \in \omega$ ;
\item $B(x_j,r) \subset B(x_{j+1},2r)~\forall j = 1,2,..,m-1$;
\item $B(x_j,2r) \Subset \Omega~\forall j=1,2,...m$ .
\end{enumerate} 
Then, thanks to (\ref{step7}), there exist $\sigma_1, \kappa _1, \sigma_m, \kappa _m$ such that
\begin{eqnarray}
&~&\int_{\Omega\cap B(x_0,R-4\epsilon)}|v(x,T)|^2dx\notag\\
&\leq& \kappa e^{\frac{\kappa }{T}}\left(\int_{B(x_0,r)} |v(x,T)|^2dx\right)^{\sigma}\left(\int_{\Omega} |v(x,0)|^2dx\right)^{1-\sigma}\notag\\
&\leq& \kappa e^{\frac{\kappa }{T}}\left(\int_{B(x_1,2r)} |v(x,T)|^2dx\right)^{\sigma}\left(\int_{\Omega} |v(x,0)|^2dx\right)^{1-\sigma}\notag\\
&\leq& \kappa e^{\frac{\kappa }{T}}\left(\kappa _1e^{\frac{\kappa _1}{T}}\left(\int_{B(x_1,r)} |v(x,T)|^2dx\right)^{\sigma_1}\left(\int_{\Omega} |v(x,0)|^2dx\right)^{1-\sigma_1}\right)^{\sigma}\left(\int_{\Omega} |v(x,0)|^2dx\right)^{1-\sigma}\notag\\
&\leq&...\notag\\
&\leq& \kappa _me^{\frac{\kappa _m}{T}}\left(\int_{B(x_m,r)} |v(x,T)|^2dx\right)^{\sigma_m}\left(\int_{\Omega} |v(x,0)|^2dx\right)^{1-\sigma_m}\notag\\
&\leq& \kappa _me^{\frac{\kappa _m}{T}}\left(\int_{\omega} |v(x,T)|^2dx\right)^{\sigma_m}\left(\int_{\Omega} |v(x,0)|^2dx\right)^{1-\sigma_m}\text{.}
\label{step8}
\end{eqnarray}
Thus, we already prove that if $\Omega \cap B(x_0,R)$ is star-shaped with respect to  $x_0$, then we obtain a local observation at one point of time, which has form (\ref{step8}). Now, in order to get the global result, we will cover $\Omega$ by a finite number of balls $B(x_0, R-4\epsilon)$ satisfying assumption that $\Omega \cap B(x_0,R)$ is star-shaped with respect to  $x_0$.\\

\textbf{Step 6} Cover $\Omega$.\\
We can see that $\Omega$ is covered by a finite number of balls $B(x_0, R-4\epsilon)$ which have one of two following properties:
\begin{enumerate}
\item $\overline{B(x_0, R)} \subset \Omega$;
\item $B(x_0, R)\cap \partial \Omega \neq \emptyset$.
\end{enumerate}
For the first case, obviously, the assumption that $\Omega \cap B(x_0,R)$ is star-shaped with respect to $x_0$ is satisfied because of the convexity of the ball $B(x_0,R)$. For the second one, we will use the result in [AEWZ] (see Theorem 8, page 2443), which says that $\Omega$ is locally star-shaped, i.e for each $\chi \in \partial \Omega$, there are $x_\chi$ in $\Omega$ and $R_\chi >0$ such that
\begin{equation*}
\chi \in B(x_\chi,R_\chi)~~ \text{and}~~ \Omega \cap B(x_\chi,R_\chi) ~~\text{is star-shaped with center} ~~x_\chi \text{.}
\end{equation*}
Thus, we can choose $x_0=x_\chi$ and $R=R_\chi$ then $\Omega \cap B(x_0,R)$ is star-shaped with respect to $x_0$.
\end{proof}

\section{ Approximate controllability at one point of time}
In \cite{LR}, Lebeau and Robbiano connect the controllability to an interpolation estimate for an elliptic system. Then, in \cite{FI}, Fursikov and Imanuvilov use a global Carleman inequality and a minimization technique to construct the control function. Recently, in \cite{FCZ}, Fern\'{a}ndez-Cara and Zuazua establish a null controllability for semilinear heat equation. In \cite{Vo}, the author succeeds in computing a control function for the cubic semilinear heat equation in a constructive way. Those results are related to the controllability in $L^2(\Omega \times (0,T))$. Here, we need to add a control at a fixed point of time, which is well studied in \cite{PWX}, but for case $p \equiv 1$. Now, our concern is approximate controllability at one point of time for linear heat equation with time-dependent coefficients. This control will lead the given data at the initial time to an origin-center ball with a small radius at some later time.\\
Denote $\mathbbm{1}_\omega$ be the characteristic function on the region $\omega$ and $\varphi(T^-)$ be the left limit of the function $\varphi$ at time $T$. Now, consider the following system.
\begin{equation}
\left\{\begin{array}{ll}
\partial_{t}\varphi-p(t)\Delta \varphi=0 & \text{in}~\Omega\times (0,2T)\setminus\{T\},\\
\varphi =0 & \text{on}~\partial \Omega \times (0,2T),\\
\varphi (\cdot,0)=\varphi ^0& \text{in}~ \Omega \text{,}\\
\varphi (\cdot,T)=\varphi (\cdot,T^{-}) + \mathbbm{1}_{\omega}{h} &\text{in} ~\Omega \text{.}
\end{array}\right.
\label{v}
\end{equation}
The next theorem consists on the existence of the control function $h$ at some fixed point of time $T$ which leads the solution at final time $2T$ getting small.
\begin{theorem}
Let $\varepsilon > 0$ and
 $\varphi ^0 \in L^2(\Omega)$,  there exists a function $h \in L^2(\omega)$ such that the solution of (\ref{v}) satisfies $\left\Vert \varphi (\cdot,2T) \right\Vert_{L^2(\Omega)}\leq \varepsilon\left\Vert \varphi ^0 \right\Vert_{L^2(\Omega)}$. Moreover, there exist constants $c_3=c_3(\Omega,\omega,p)>0$ and $c_4=c_4(\Omega,\omega,p)>0$ such that
\begin{equation} \left\Vert h\right\Vert_{L^2(\omega)} \leq \frac {c_3e^{\frac{c_3}{T}}}{\varepsilon^{c_4}}\left\Vert \varphi ^0 \right\Vert_{L^2(\Omega)}\text{.}\label{confun}\end{equation}
\label{control}
\end{theorem}
 \textbf{Proof of Theorem \ref{control}}
  \begin{proof}
\textbf{Step 1:} Define a functional which has a unique minimizer.\\
 Let $c_1$ and $c_2$ be the constants from Corollary \ref{colyoung}, put $k:=\sqrt{c_1e^{\frac{c_1}{T}}\frac{1}{\varepsilon^{2c_2}}}$. We consider the following functional
 \begin{equation} J(\mathfrak{u}_0)=\frac{k^2}{2}\Vert \mathfrak{u}(\cdot,T)\Vert^2_{L^2(\omega)}+\frac{\varepsilon^2}{2}\Vert \mathfrak{u}_0\Vert^2_{L^2(\Omega)}-\int_{\Omega} \varphi^0(x)\mathfrak{u}(x,2T)dx\end{equation}
 where $\mathfrak{u}(x,t)$ is the solution of 
 \begin{equation}
\left\{\begin{array}{ll}
\partial_{t}\mathfrak{u}-p(t)\Delta \mathfrak{u}=0 & \text{in}~\Omega\times (0,T) \text{,}\\
\mathfrak{u}=0 & \text{on}~\partial \Omega \times (0,T) \text{,}\\
\mathfrak{u}(\cdot,0)=\mathfrak{u}_0 &\text{in}~ \Omega \text{.}
\label{J1}	
\end{array}\right.
\end{equation}
 Notice that $J$ is a strictly convex, $C^1$ and coercive. Therefore, $J$ has a unique minimizer $\Phi_0 \in L^2(\Omega)$ such that $J(\Phi_0)=\min\limits_{\mathfrak{u}_0\in L^2(\Omega)} J(\mathfrak{u}_0)$. It implies that $J'(\Phi_0)\mathfrak{z}_0 = 0$ for any $\mathfrak{z}_0 \in L^2(\Omega)$, i.e the following estimate holds for any $\mathfrak{z}_0$:
  \begin{equation}
 k^2\Vert \Phi(\cdot,T)\Vert_{L^2(\omega)}\Vert \mathfrak{z}(\cdot,T)\Vert_{L^2(\omega)}+\varepsilon^2\Vert \Phi_0\Vert_{L^2(\Omega)}\Vert \mathfrak{z}_0\Vert_{L^2(\Omega)}-\int_{\Omega} \varphi^0(x)\mathfrak{z}(x,2T)dx =0
 \label{J'}
 \end{equation}
 where $\Phi(x,t)$ and $\mathfrak{z}(x,t)$ are respectively the solution of (\ref{J1}) corresponding to $\Phi_0:=\Phi(\cdot,0)$ and $\mathfrak{z}_0:=\mathfrak{z}(\cdot,0)$. \\
 
\textbf{Step 2:} Construct a control function.\\
 On the other hand, multiplying $\partial_{t}\varphi -p(t)\Delta \varphi =0$ by $\Phi(\cdot,2T-t)$ and integrating over $\Omega$, we get
 \begin{equation}
 \int_\Omega \partial_t \varphi (x,t) \Phi(x,2T-t)dx - p(t) \int_\Omega \Delta \varphi (x,t) \Phi(x,2T-t)dx =0
 \text{.}
 \label{114}\end{equation}
Integrating by parts (\ref{114}) two times with the fact that $\varphi = \Phi =0$ on $\partial \Omega$, one has
  \begin{equation}
 \int_\Omega \partial_t \varphi (x,t) \Phi(x,2T-t)dx - p(t)\int_\Omega \varphi (x,t) \Delta \Phi(x,2T-t)dx =0
 \text{.}\end{equation}
 Since $\partial_t \Phi -p(t) \Delta \Phi =0$, it yields
  \begin{equation}
 \int_\Omega \partial_t \varphi (x,t) \Phi(x,2T-t)dx - \int_\Omega \varphi (x,t) \partial_t \Phi(x,2T-t)dx =0
 \text{.}\end{equation}
 It implies that
 \begin{equation}
 \frac{d}{dt}\int_\Omega \varphi (x,t)\Phi(x,2T-t)dx =0
 \label{dt}
 \text{.}\end{equation}
 Now, by integrating (\ref{dt}) over $(0,T)$, we obtain
 \begin{equation}
 \int_\Omega \varphi (x,0)\Phi(x,2T)dx =  \int_\Omega \varphi (x,T^-)\Phi(x,T)dx
 \label{11}
 \text{.}\end{equation}
 Integrating again (\ref{dt}) over $(T,2T)$ forces
 \begin{eqnarray}
 \int_\Omega \varphi (x,2T)\Phi(x,0)dx&=& \int_\Omega \varphi (x,T)\Phi(x,T)dx\notag\\
 &=& \int_\Omega \varphi (x,T^-)\Phi(x,T)dx + \int_\Omega h(x)\Phi(x,T)dx
 \text{.}\end{eqnarray}
 Combining the above equality with (\ref{11}), we conclude that
 \begin{equation}
 \int_\Omega \varphi (x,2T)\Phi_0dx =  \int_\Omega \varphi ^0(x)\Phi(x,2T)dx + \int_\Omega h(x)\Phi(x,T)dx
 \end{equation}
 that is
  \begin{equation}
 -\int_\Omega h(x)\Phi(x,T)dx+\int_\Omega \varphi (x,2T)\Phi_0dx -\int_\Omega \varphi ^0(x)\Phi(x,2T)dx=0 \label{22}\text{.}
 \end{equation}
In addition, by choosing $\mathfrak{z}_0 = \Phi_0$ in (\ref{J'}), it follows that
 \begin{equation}
 k^2\Vert \Phi(\cdot,T)\Vert^2_{L^2(\omega)}+\varepsilon^2\Vert \Phi_0\Vert^2_{L^2(\Omega)}-\int_{\Omega} \varphi ^0(x)\Phi(x,2T)dx =0 \text{.}
 \label{23}
 \end{equation}
Thus from (\ref{22}) and (\ref{23}), if we choose $h(x)=-k^2\Phi(x, T)$ then $\varphi (x,2T) = \varepsilon^2 \Phi_0(x)$. Moreover, using the Cauchy-Schwarz inequality
 \begin{eqnarray}
 \int_{\Omega} |v^0(x)\Phi(x,2T)|dx \leq \left(\int_{\Omega} |\varphi ^0(x)|^2dx\right)^{\frac{1}{2}}\left(\int_{\Omega} |\Phi(x,2T)|^2dx\right)^{\frac{1}{2}}
 \text{,}\end{eqnarray}
 we get
  \begin{eqnarray}
 k^2\Vert \Phi(\cdot,T)\Vert^2_{L^2(\omega)}+\varepsilon^2\Vert \Phi_0\Vert^2_{L^2(\Omega)}\leq \Vert \varphi ^0\Vert_{L^2(\Omega)}\Vert \Phi(\cdot,2T)\Vert_{L^2(\Omega)}
  \text{.}
  \label{tamtam}
  \end{eqnarray}
 On the other hand, Corollary \ref{colyoung} gives
 \begin{equation} \Vert \Phi(\cdot,T)\Vert^2_{L^2(\Omega)} \leq  k^2\Vert \Phi(\cdot,T)\Vert^2_{L^2(\omega)}+\varepsilon^2\Vert \Phi(\cdot,0)\Vert^2_{L^2(\Omega)}\text{.}\end{equation}
Furthermore, the fact that $\Vert \Phi(\cdot,2T)\Vert^2_{L^2(\Omega)}\leq \Vert \Phi(\cdot,T)\Vert^2_{L^2(\Omega)}$ provides us
 \begin{equation} \Vert \Phi(\cdot,2T)\Vert^2_{L^2(\Omega)} \leq  k^2\Vert \Phi(\cdot,T)\Vert^2_{L^2(\omega)}+\varepsilon^2\Vert \Phi(\cdot,0)\Vert^2_{L^2(\Omega)}\text{.}
 \label{k}
 \end{equation}
 Thus, combining (\ref{tamtam}) and (\ref{k}), we conclude that
 \begin{eqnarray*}
 &~&k^2\Vert \Phi(\cdot,T)\Vert^2_{L^2(\omega)}+\varepsilon^2\Vert \Phi_0\Vert^2_{L^2(\Omega)}\leq \Vert \varphi ^0\Vert_{L^2(\Omega)}\left(k^2\Vert \Phi(\cdot,T)\Vert^2_{L^2(\omega)}+\varepsilon^2\Vert \Phi_0\Vert^2_{L^2(\Omega)}\right)^{\frac{1}{2}}
 \text{.}\end{eqnarray*}
It implies that
\begin{equation}k^2\Vert \Phi(\cdot,T)\Vert^2_{L^2(\omega)}+\varepsilon^2\Vert \Phi_0\Vert^2_{L^2(\Omega)} \leq \Vert \varphi ^0\Vert^2_{L^2(\Omega)}\text{.}\end{equation}
This is equivalent to
\begin{equation}
\frac{1}{k^2}\Vert h\Vert^2_{L^2(\omega)} + \frac{1}{\varepsilon^2}\Vert \varphi (\cdot,2T)\Vert^2_{L^2(\Omega)} \leq \Vert \varphi ^0\Vert^2_{L^2(\Omega)}\text{.}\end{equation}
Thus, we get the desired estimate (\ref{confun}) with $c_3:=\max\{\sqrt{c_1},\frac{c_1}2\}$ and $c_4:=c_2$. This completes the proof.
\end{proof}
\section{The local backward - Proof of the Theorem \ref{local}}
For the case $\omega \Subset \Omega$, we need to use the controllability result at one point of time (Theorem \ref{control}) in order to get the information of the solution on the whole domain from the known data on the subdomain. In detail, the proof of Theorem \ref{local} is structured as: Step 1 will provide us the approximate data of $u(\cdot, 3T)$ on whole domain dues to the controllability result; Then, by some computation technique, we make appear $f$ in Step 2; Lastly, applying the global backward result in Theorem \ref{global}, Step 3 will complete the proof with the construction of the initial data.
\begin{proof}
\textbf{Step 1:} Use controllability result to link the knowledge on the whole domain $\Omega$ and on the subdomain $\omega$.\\
Now, for each $i=1,2,...$, Theorem \ref{control} says that for any $\varepsilon >0$, there exists $h_i \in L^2(\omega)$ such that the solution of
\begin{equation}
\left\{\begin{array}{ll}
\partial_{t}\varphi _i-p(t)\Delta \varphi _i=0 & \text{in}~\Omega\times (0,2T)\setminus\{T\},\\
\varphi _i=0 & \text{on}~\partial \Omega \times (0,2T),\\
\varphi _i(\cdot,0)=e_i& \text{in}~ \Omega \text{,}\\
\varphi _i(\cdot,T)=\varphi _i(\cdot,T^{-}) + \mathbbm{1}_{\omega}{h_i} &\text{in} ~\Omega
\end{array}\right.
\label{vi}
\end{equation}
satisfies $\left\Vert \varphi _i(\cdot,2T) \right\Vert_{L^2(\Omega)}\leq \varepsilon$ for any $i\geq 1$. Recall that $e_i (i=1,2,...)$ is the eigenfunction of Laplace operator. Moreover, there exist constants $c_3, c_4 > 0$ such that the following estimate holds 
\begin{equation}
\left\Vert h_i\right\Vert_{L^2(\omega)} \leq \frac{c_3e^{\frac{c_3}{T}}}{\varepsilon^{c_4}}~~~\forall i \geq 1\text{.}
\label{fi}
\end{equation}
Multiplying both sides of the equation $\partial_{t}\varphi _i-p(t)\Delta \varphi _i=0$ by $u(\cdot, 2T-t)$ and integrating over $\Omega$, one gets
\begin{equation} \frac{d}{dt}\int_{\Omega} \varphi _i(x,t) u(x,2T-t) dx = 0\text{.}\label{222}\end{equation}
Integrating (\ref{222}) over $(0,T)$ and $(T,2T)$ respectively and using the fact that $\varphi _i(\cdot,T)=\varphi _i(\cdot,T^{-}) + \mathbbm{1}_{\omega}{h_i}$, one has
\begin{equation}\int_\Omega \varphi _i(x,2T)u(x,0)dx = \int_\Omega \varphi _i(x,0) u(x,2T)dx + \int_\omega h_i(x) u(x,T)dx\text{.}\label{ice}\end{equation}
Replacing $$u(\cdot,2T) = \sum\limits_{j=1}^{\infty}e^{-\lambda_j \int_0^{2T} p(s)ds}\int_\Omega u(x,0)e_j(x)dx e_j$$ and $\varphi _i(\cdot,0) = e_i$ in (\ref{ice}), we get
\begin{equation}\int_\Omega \varphi _i(x,2T)u(x,0)dx=e^{-\lambda_i \int_0^{2T} p(s)ds}\left(\int_\Omega u(x,0)e_i(x)dx\right) + \int_\omega h_i(x) u(x,T)dx\label{coucou}\text{.}\end{equation}
Now, multiplying both sides of (\ref{coucou}) by $e^{-\lambda_i \int_{2T}^{3T} p(s)ds}e_i$ and take the sum from $i=1$ to $\infty$, one has
\begin{eqnarray}
&~&u(\cdot,3T)+ \sum\limits_{i=1}^{\infty} e^{-\lambda_i \int_ {2T}^{3T}p(s)ds}\int_\omega h_i(x) u(x,T)dx e_i\notag\\
&=&\sum\limits_{i=1}^{\infty} e^{-\lambda_i \int_{2T}^{3T} p(s)ds} \int_\Omega \varphi_i(x,2T)u(x,0)dx e_i\text{.}
\end{eqnarray}
It implies from $\left\Vert \varphi_i(\cdot,2T) \right\Vert_{L^2(\Omega)}\leq \varepsilon$ that
\begin{eqnarray}
&~&\left\Vert u(\cdot,3T) + \sum\limits_{i=1}^{\infty} e^{-\lambda_i \int_{2T}^{3T} p(s)ds}\int_\omega h_i(x) u(x,T)dx e_i\right\Vert_{L^2(\Omega)}\notag\\
& \leq &\varepsilon \left\Vert u(\cdot,0)\right\Vert_{L^2(\Omega)} \left( \sum\limits_{i=1}^{\infty} e^{-2\lambda_i \int_{2T}^{3T} p(s)ds}\right)^{\frac{1}{2}}
\text{.}
\label{3T}
\end{eqnarray}

\textbf{Step 2:} Make appear $f$.\\
Now, we will make appear $f$ by using a triangle inequality
\begin{eqnarray}
&~&\left\Vert u(\cdot, 3T) + \sum\limits_{i=1}^{\infty} e^{-\lambda_i \int_{2T}^{3T} p(s)ds}\int_\omega h_i(x)f(x)dx e_i\right\Vert_{L^2(\Omega)}\notag\\
&\leq&  \left\Vert u(\cdot, 3T) + \sum\limits_{i=1}^{\infty} e^{-\lambda_i \int_{2T}^{3T} p(s)ds}\int_\omega h_i(x) u(x,T)dx e_i \right\Vert_{L^2(\Omega)} \notag\\
&~&+ \left\Vert \sum\limits_{i=1}^{\infty} e^{-\lambda_i \int_{2T}^{3T} p(s)ds}\int_\omega h_i(x) \left(u(x,T)-f\right)dx e_i \right\Vert_{L^2(\Omega)}\text{.}
\end{eqnarray}
We got the first estimate in Step 1, let us estimate the second one.
\begin{eqnarray}
&~&\left\Vert \sum\limits_{i=1}^{\infty} e^{-\lambda_i \int_{2T}^{3T} p(s)ds}\int_\omega h_i(x) \left(u(x,T)-f\right)dx e_i \right\Vert_{L^2(\Omega)}\notag\\
&\leq& \left(\sum\limits_{i=1}^{\infty} e^{-2\lambda_i \int_{2T}^{3T} p(s)ds}\right)^{\frac{1}{2}} \left\Vert h_i\right\Vert_{L^2(\omega)}\left\Vert u(\cdot,T)-f\right\Vert_{L^2(\omega)}\notag\\
&\leq& \left(\sum\limits_{i=1}^{\infty} e^{-2\lambda_i \int_{2T}^{3T} p(s)ds}\right)^{\frac{1}{2}} \frac{c_3e^{\frac{c_3}{T}}}{\varepsilon^{c_4}}\delta \text{.}
\label{hehehe}
\end{eqnarray}
The last inequality in (\ref{hehehe}) is followed from (\ref{fi}) and the assumption (\ref{globalz}). 
Thus, from (\ref{3T}) and (\ref{hehehe}), we can conclude that
\begin{eqnarray}
&~&\left\Vert u(\cdot, 3T) + \sum\limits_{i=1}^{\infty} e^{-\lambda_i \int_{2T}^{3T} p(s)ds}\int_\omega h_i(x)f(x)dx e_i\right\Vert_{L^2(\Omega)}\notag\\
&\leq&  \left( \sum\limits_{i=1}^{\infty} e^{-2\lambda_i \int_{2T}^{3T} p(s)ds}\right)^{\frac{1}{2}}\left[\varepsilon \left\Vert u(\cdot,0)\right\Vert_{L^2(\Omega)}+ \frac{c_3e^{\frac{c_3}{T}}}{\varepsilon^{c_4}}\delta\right]\text{.}
\label{cat}
\end{eqnarray}
It is known that the function $x \mapsto ax+bx^{-s}~~(a,b,s>0)$ gets minimum at $x_0 = \left(\frac{sb}{a}\right)^{\frac 1{s+1}}$. Hence, in order to minimize the right-hand side of (\ref{cat}) we choose $\varepsilon$ such that
\begin{equation} \varepsilon = \left(\frac{c_4c_3e^{\frac{c_3}{T}}}{\left\Vert u(\cdot,0) \right\Vert_{L^2(\Omega)}}\delta\right)^{\frac{1}{1+c_4}}\text{.}\end{equation}
Therefore, (\ref{cat}) becomes
\begin{eqnarray}
&~&\left\Vert u(\cdot, 3T) + \sum\limits_{i=1}^{\infty} e^{-\lambda_i \int_{2T}^{3T} p(s)ds}\int_\omega h_i(x)f(x)dx e_i \right\Vert_{L^2(\Omega)} \notag\\
&\leq & \left( \sum\limits_{i=1}^{\infty} e^{-2\lambda_i \int_{2T}^{3T} p(s)ds}\right)^{\frac{1}{2}}\left(\frac{c_4c_3e^{\frac{c_3}{T}}\delta}{\left\Vert u(\cdot,0) \right\Vert_{L^2(\Omega)}}\right)^{\frac{1}{1+c_4}}\left(1+\frac 1{c_4}\right)\left\Vert u(\cdot,0)\right\Vert_{L^2(\Omega)} \notag\\
&\leq & \left( \sum\limits_{i=1}^{\infty} e^{-2\lambda_i p_2T}\right)^{\frac{1}{2}}\left(\frac{c_4c_3e^{\frac{c_3}{T}}\delta}{\left\Vert u(\cdot,0) \right\Vert_{L^2(\Omega)}}\right)^{\frac{1}{1+c_4}}\left(1+\frac 1{c_4}\right)\left\Vert u(\cdot,0)\right\Vert_{L^2(\Omega)} \text{.}
\label{sat}
\end{eqnarray}
Using the fact that $e^{-x} \leq \left(\frac \gamma x\right)^{\gamma}~\forall x>0~\forall \gamma >0$ and $\lambda_i \approx i^{\frac 2n}~~\forall i\geq 1$ ( by the Weyl formula), we get: 
\begin{eqnarray*}
\sum\limits_{i=1}^\infty e^{-2\lambda_ip_2T} &\leq& \left(\frac \gamma {2p_2T}\right)^\gamma \sum\limits_{i=1}^\infty \frac 1{\lambda_i^\gamma}\leq \left(\frac \gamma {2p_2T}\right)^\gamma \sum\limits_{i=1}^\infty \frac 1{\left(i^{\frac 2 n}\right)^\gamma}\text{.}
\end{eqnarray*}
When $\gamma > \frac n2$, there exists a constant $S >0$ such that $\sum\limits_{i=1}^\infty \frac 1{\left(i^{\frac 2 n}\right)^\gamma} = S$. Hence
\begin{eqnarray*}
\sum\limits_{i=1}^\infty e^{-2\lambda_ip_2T} \leq e^{\frac {\gamma^2} {2p_2T} }S\text{.}
\end{eqnarray*}
One can conclude that (\ref{sat}) can be written as below
\begin{eqnarray}
\left\Vert u(\cdot, 3T) + \sum\limits_{i=1}^{\infty} e^{-\lambda_i \int_{2T}^{3T} p(s)ds}\int_\omega h_i(x)f(x)dx e_i \right\Vert_{L^2(\Omega)}\leq K_1e^{\frac {K_1}T}\Vert u(\cdot,0)\Vert_{L^2(\Omega)}^{1-k_1}\delta^{k_1}
\end{eqnarray}
for some positive constants $K_1=K_1(\Omega, \omega,p)>1$ and $k_1=k_1(\Omega, \omega, p)\in (0,1)$.\\

\textbf{Step 3:} Apply the global backward result.\\
From the facts that $h_i\in L^2(\omega)$, $z\in L^2(\omega)$ and $\sum\limits_{i=1}^{\infty} e^{-\lambda_i \int_{2T}^{3T} p(s)ds}<\infty$, one gets
\begin{equation} - \sum\limits_{i=1}^{\infty} e^{-\lambda_i \int_{2T}^{3T} p(s)ds}\int_\omega h_i(x) f(x)dx e_i \in L^2(\Omega) 
\text{.}\end{equation}
Thus, Theorem \ref{global} gives us the following estimate, where $\delta$ in (\ref{globalz}) is replaced by $K_1e^{\frac {K_1}T}\Vert u(\cdot,0)\Vert_{L^2(\Omega)}^{1-k_1}\delta^{k_1}$
\begin{equation}
\left\Vert u(\cdot,0) - g \right\Vert _{L^2(\Omega)}  \leq \frac{\sqrt{(1+\zeta)p_2T}\left\Vert u(\cdot,0)\right\Vert_{H_0^1(\Omega)}}{\sqrt{\ln \left(\sqrt{2\zeta\lambda_1 p_2T}\frac{\left\Vert u(\cdot,0)\right\Vert_{L^2(\Omega)}}{K_1e^{\frac {K_1}T}\Vert u(\cdot,0)\Vert_{L^2(\Omega)}^{1-k_1}\delta^{k_1}}\right)}}\text{.}
\end{equation}
for any $ \zeta > \frac{K^2_1e^{\frac {2K_1}T}}{2\lambda_1p_2T}\left(\frac{\delta}{\Vert u(\cdot,0)\Vert^2_{L^2(\Omega)}}\right)^{2k_1} $. 
With the assumption that $\delta<\Vert u(\cdot,0)\Vert_{L^2(\Omega)}$, we can choose $\zeta$ as below
\begin{equation}
\zeta = \frac {K_1^2e^{2\frac{K_1}{T}}}{2\lambda_1p_2T}
\end{equation}
in order to get
\begin{equation}
\left\Vert u(\cdot,0) - g \right\Vert _{L^2(\Omega)}  \leq \frac{\left(1+ \frac {K_1e^{\frac{K_1}{T}}}{\sqrt{2\lambda_1p_2T}}\right)\sqrt{p_2T}\left\Vert u(\cdot,0)\right\Vert_{H_0^1(\Omega)}}{\sqrt{k_1}\sqrt{\ln \frac {\Vert u(\cdot,0)\Vert_{L^2(\Omega)}}{\delta}}}\text{.}
\end{equation}
Using the fact that $x\leq e^x~~\forall x>0$, there exists a constant $C>0$ such that
\begin{equation}
\left\Vert u(\cdot,0) - g \right\Vert _{L^2(\Omega)}  \leq \frac{Ce^{\frac{C}{T}}\sqrt{T}\left\Vert u(\cdot,0)\right\Vert_{H_0^1(\Omega)}}{\sqrt{\ln \frac {\Vert u(\cdot,0)\Vert_{L^2(\Omega)}}{\delta}}}\text{.}
\end{equation}
This completes the proof.
\end{proof}
\subsection{Appendix}
\textbf{Proof of (\ref{appendix}) in Remark \ref{remark2}}\\
Let us remind the observation estimate (\ref{holder}) in Theorem \ref{nonconvex} (see Section 3)
\begin{equation}
 \Vert u(\cdot,T)\Vert_{L^2(\Omega)} \leq K_2e^{\frac{K_2}{T}} \Vert u(\cdot,0)\Vert_{L^2(\Omega)}^{1-k_2}\Vert u(\cdot,T)\Vert_{L^2(\omega)}^{k_2}
 \end{equation}
for some $K_2=K_2(\Omega, \omega, p)>0$ and $k_2=k_2(\Omega, \omega, p) \in (0,1)$. Combining to the known backward estimate below
\begin{equation}
\left\Vert u(\cdot,0)\right\Vert _{L^2(\Omega)}  \leq  e^{\frac{p_2T \left\Vert u(\cdot,0)\right\Vert^2_{H_0^1(\Omega)}}{\left\Vert u(\cdot,0)\right\Vert^2_{L^2(\Omega)}}} \left\Vert u(\cdot,T)\right\Vert _{L^2(\Omega)}\text{,}
\end{equation}
we get
\begin{equation}
\left\Vert u(\cdot,0)\right\Vert _{L^2(\Omega)}  \leq  e^{\frac{p_2T \left\Vert u(\cdot,0)\right\Vert^2_{H_0^1(\Omega)}}{\left\Vert u(\cdot,0)\right\Vert^2_{L^2(\Omega)}}}K_2e^{\frac{K_2}{T}} \Vert u(\cdot,0)\Vert_{L^2(\Omega)}^{1-k_2}\Vert u(\cdot,T)\Vert_{L^2(\omega)}^{k_2}\text{.}
 \end{equation}
 It is equivalent to 
 \begin{equation}
\left(\frac{\left\Vert u(\cdot,0)\right\Vert _{L^2(\Omega)}}{\left\Vert u(\cdot,T)\right\Vert _{L^2(\omega)}}\right)^{k_2} \leq e^{\frac{p_2T \left\Vert u(\cdot,0)\right\Vert^2_{H_0^1(\Omega)}}{\left\Vert u(\cdot,0)\right\Vert^2_{L^2(\Omega)}}}K_2e^{\frac{K_2}T}\text{.}
\end{equation}
It follows that
\begin{equation}
\frac{\left\Vert u(\cdot,0)\right\Vert _{L^2(\Omega)}}{\left\Vert u(\cdot,T)\right\Vert _{L^2(\omega)}} \leq e^{\left(\frac{p_2}{k_2}\frac{T \left\Vert u(\cdot,0)\right\Vert^2_{H_0^1(\Omega)}}{\left\Vert u(\cdot,0)\right\Vert^2_{L^2(\Omega)}}+\frac{K_2}{k_2}\frac 1T+\frac 1 {k_2} \ln K_2\right)}\text{.}
\end{equation}
Using the fact that $\frac{\left\Vert u(\cdot,0)\right\Vert^2_{H_0^1(\Omega)}}{\left\Vert u(\cdot,0)\right\Vert^2_{L^2(\Omega)}}\geq \lambda_1$, there exists a constant $C=\sqrt{\max\{\frac{p_2}{k_2}, \frac {K_2}{k_2\lambda_1}\}}>0$ such that
\begin{equation}
\frac{\left\Vert u(\cdot,0)\right\Vert _{L^2(\Omega)}}{\left\Vert u(\cdot,T)\right\Vert _{L^2(\omega)}} \leq e^{C^2\left(1+T+\frac 1{T}\right) \frac{\left\Vert u(\cdot,0)\right\Vert^2_{H_0^1(\Omega)}}{\left\Vert u(\cdot,0)\right\Vert^2_{L^2(\Omega)}}}\text{.}
\end{equation}
Thus
\begin{equation}
\left\Vert u(\cdot,0)\right\Vert _{L^2(\Omega)}  \leq \frac{C\sqrt{1+T+\frac{1}{T}}\Vert u(\cdot, 0)\Vert_{H_0^1(\Omega)}}{\sqrt{\ln \frac{\Vert u(\cdot,0)\Vert_{L^2(\Omega)}}{\Vert u(\cdot,T)\Vert_{L^2(\omega)}}}}\text{.}
\end{equation}

 \end{document}